\documentclass[11 pt,Rectoverso]{article}  
\font\dsrom=dsrom10 scaled 1400

\def\di{/\!\!/}

\def\proof{\noindent{\bf Proof. }}

\def \1{\textrm{\dsrom{1}}}
\def\wh{\widehat}

\usepackage[latin1]{inputenc}
\usepackage{graphicx,psfrag,amsfonts,bbm}
\usepackage{amsmath,indentfirst,qsymbols}

\def \app#1#2#3#4#5{\begin{array}{rccl} #1:&#2&\longrightarrow&#3\\ &#4&\longmapsto&#5\end{array}}
\def \sk{{\bf s}^{(\kappa)}}
\def \be{\begin{eqnarray*}}
\def \ee{\end{eqnarray*}}
\def \ben{\begin{eqnarray}}
\def \een{\end{eqnarray}}
\def \bq{\begin{equation}}
\def \eq{\end{equation}}

\def \build#1#2#3{\mathrel{\mathop{\kern 0pt#1}\limits_{#2}^{#3}}}
\def \cro#1{\llbracket#1\rrbracket}
\def \croc#1{\rrbracket#1\llbracket}

\def \floor#1{\lfloor#1\rfloor}

\def \eref#1{(\ref{#1})}
\def \sous#1#2{\mathrel{\mathop{\kern 0pt#1}\limits_{#2}}}
\def \sur#1#2{\mathrel{\mathop{\kern 0pt#1}\limits^{#2}}}

\def \dd{\xrightarrow[n]{(d)}}
\def \proba{\xrightarrow[n]{proba.}}
\def \as{\xrightarrow[n]{a.s.}}
\def \ov#1{\overrightarrow{#1}}

\def \bh{{\bf h}}

\def \br{{\bf r}}

\def \bG{{\bf G}}
\def \bGn{{\bf G}^{(n)}}

\def \bg{{\bf g}}
\def \sa{{\sf a}}
\def \sx{{\sf x}}
\def \se{{\sf e}}

\def \captionn#1{\begin{center}\begin{minipage}{14cm}\sf\caption{\small #1}\end{minipage}\end{center}}

\def \T{{\mathcal T}}

\def \hh#1{^{(#1)}}

\def \dis{\displaystyle}
\def \tend{\longrightarrow}

\def \b{\big}
\def \bar{\overline}

\def \precc{\preccurlyeq}
\def \m{{\cal M}}
\def \l{\left}
\def \r{\right}

\def \B{{\cal B}}
\def \root{\varnothing}
\def \ni{\mathbb{N}^I}
\newcommand{\bt}{{\bf T}}

\font\dsrom=dsrom10 scaled 1400
\newqsymbol{`I}{\textrm{\dsrom{1}}}
\newqsymbol{`P}{\mathbb{P}}
\newqsymbol{`Q}{\mathbb{Q}}
\newqsymbol{`R}{\mathbb{R}}
\newqsymbol{`Z}{\mathbb{Z}}
\newqsymbol{`E}{\mathbb{E}}
\newqsymbol{`e}{\varepsilon}
\newqsymbol{`a}{\alpha}
\newqsymbol{`t}{\tau}
\newqsymbol{`o}{\omega}
\newqsymbol{`O}{\Omega}
\newqsymbol{`N}{\mathbb{N}}

\graphicspath{{../DESSINS/}}

\DeclareMathOperator{\var}{var}

\DeclareMathOperator{\supp}{supp}

\DeclareMathOperator{\cov}{cov}
\DeclareMathOperator\Sub{{\cal S}}
\DeclareMathOperator\fa{fa}

\begin{document}
\newtheorem{lem}{Lemma}
\newtheorem{defi}{Definition}
\newtheorem{pro}[lem]{Proposition}
\newtheorem{theo}[lem]{Theorem}
\newtheorem{cor}[lem]{Corollary}
\newtheorem{remi}{Remark\rm}{\rm}
\newtheorem{com}{Comments\rm}{\rm}
\newtheorem{exe}{Examples \rm}{\rm}
\newenvironment{paren}%
{\footnotesize}%

\begin{center}\LARGE{\bf Globally centered discrete snakes}\\\medskip\normalsize
 \textrm{\Large Jean-Fran\c{c}ois Marckert}\\
 \textrm{CNRS, LaBRI}\\
\textrm{Universit\'e Bordeaux 1}\\
   \textrm{351 cours de la Libération}\\
 \textrm{33405 Talence cedex, FRANCE}\\
 \end{center}
\begin{abstract}
We consider branching random walks built on Galton-Watson trees with offspring distribution having a bounded support, conditioned to have $n$ nodes, and their rescaled convergences  to the Brownian snake. We exhibit a notion of ``globally centered discrete snake'' that extends the usual settings in which the displacements are supposed centered. We show that under some additional moment conditions, when $n$ goes to $+\infty$, ``globally centered discrete snakes'' converge to the Brownian snake. 
The proof relies on a precise study of the ``lineage'' of the nodes in a Galton-Watson tree conditioned by the size, and their links with a multinomial process. Some consequences concerning Galton-Watson trees conditioned by the size are also derived.
\end{abstract}
\medskip 
\normalsize

\section{Introduction}

\subsection{A model of centered discrete snake}
We first begin with the formal description of the notion of trees and branching random walks. \par
Let $T_{\infty}=\{\varnothing\}\cup\bigcup_{n\geq 1}`N^\star{}^n$ be the set of finite words on the alphabet $\mathbb{N}^{\star}=\{1,2,\dots\}$.  For $u=u_1\ldots u_n$, and $v=v_1\ldots v_m\in T_{\infty}$, we let $uv=u_1\ldots u_nv_1\ldots v_m$ be the concatenation of the words $u$ and $v$ (by convention $\varnothing u=u\varnothing=u$). 
Following Neveu, we call planar tree $T$ a subset of 
$T_{\infty}$ containing the root $\varnothing$, and such that if 
$ui\in T$, then $u\in T$ and for all $j\in\cro{1,i}$, $uj\in T$. 
The elements of a tree are called nodes or vertices. For $i\neq j$, the nodes $ui$ and $uj$ are called brothers and $u$ their father. We let $c_u(T)=\max\{i:ui\in T\}$ be the number of children of $u$. A node without any child is called a leaf, and we denote by $\partial T$ the set of leaves of $T$.   
If $v \not= \varnothing$, we say that $uv$ is a descendant of $u$ and $u$ is an ancestor of $uv$. An edge is a pair $\{u,v\}$  where $u$ is the father of $v$.    
A path $\cro{u,v}$  between the nodes $u$ and $v$ in a tree $T$ is the (minimal) sequence of nodes $u:=u_0,\dots,u_j:=v$ such that for any $i\in\cro{0,j-1}$, $\{u_i,u_{i+1}\}$ is an edge. Set also $\croc{u,v}=\cro{u,v}\setminus\{u,v\}$ and similar notation for $\llbracket u,v \llbracket$ and for $\rrbracket u,v \rrbracket$.
The distance $d_T$, or simply $d$, is the usual graph distance. The depth of $u$ is $|u|=d(\varnothing,u)$. The cardinality of $T$ is denoted by $|T|$, and we let $\T$ (resp. $\T_n$) be the set of planar trees (resp. with $n$ edges, i.e. $n+1$ vertices). \par

A \it branching walk \rm is a pair $(T,\ell)$ where $T$ is a tree called the underlying tree and $\ell$, the label function, is an application from $T$ taking its values in $`R$. In other words it is a tree in which every vertex owns a real label. We let $\B$ be the set of branching walks, and $\B_n$ be the branching walks associated with trees from $\T_n$.\par

We introduce now some randomness and construct a probability distribution on $\B$ and on $\B_n$. The set of underlying trees is endowed with the distribution of the family tree of a Galton-Watson (GW) process with offspring distribution ${\bf \mu}=(\mu_k)_{k\geq 0}$ starting from one individual. We denote by  $\bt$ a random tree under this distribution.
The distribution of the labels is defined as follows. 
Consider $(\nu_{k})_{k\in\{1,2,\dots\}}$ a family of distributions, where $\nu_k$ is a distribution on $`R^k$. The labels are defined conditionally on the underlying tree $\bt$~: Set $\ell(\root)=0$, and for any $u\in \bt\setminus \partial \bt$, consider \[X_u:=\b(\ell(u1)-\ell(u),\dots,\ell(uc_u(\bt))-\ell(u)\b),\]
the evolution-vector of the labels between $u$ and its children. 
Conditionally on $\bt$, we assume that the r.v.  $X_u$ are independent, and that $X_u$ has distribution $\nu_{c_u(\bt)}$. This determines a distribution on $\B$, denoted by $`P$.
For example, if $\nu_k$ is the uniform distribution on $\{-1,+1\}^k$ for any $k>0$, then the r.v.  $\ell(u1)-\ell(u),\dots,\ell(uc_u(\bt))-\ell(u)$ are independent with common distribution $\frac12(\delta_{+1}+\delta_{-1})$ ($\delta_x$ stands for the Dirac mass at $x$). In the case where $\nu_k$ is the uniform distribution on $\{(1,\dots,k),(-1,\dots,-k)\}$, the r.v.  $\ell(ui)-\ell(u)$ and $\ell(uj)-\ell(u)$ are not independent and do not have the same distribution.\par

We define now two sets of assumptions $({\rm H}_1)$ and $({\rm H}_2)$  that will be assumed to be satisfied in most of our results.
 $({\rm H}_1)$ is the conditions that $\mu$ is non-degenerate critical and has a bounded support~:
\[({\rm H}_1):=\Big(
\mu_0+\mu_1\neq 1,
\sum_{k\geq 0} k\mu_k=1,
\textrm{there exists }K>0 \textrm{ s.t. } \sum_{k\leq K}\mu_k=1
\Big).\]
Under $({\rm H}_1)$ the variance $\sigma^2_{\mu}$ of $\mu$ is finite and non zero. 
The bounded support condition is quite a strong restriction but considering non-bounded distribution leads to non-trivial complications, and we were unable to extend to that case the most important results.

Let  $Y^{(k)}=(Y_{k,1},\dots,Y_{k,k})$ be $\nu_k$-distributed, and let $m_{k,j}$ and $\sigma_{k,j}^2$ be the mean and the variance of $Y_{k,j}$. We call \it global mean \rm and \it global variance \rm  of the branching random walk,
\[{\bf m}=\sum_{k\geq 1}\sum_{j=1}^k \mu_k m_{k,j},~~ \textrm{ and }~~ {\bf \beta}^2=\sum_{k\geq 1}\sum_{j=1}^k \mu_k  `E(Y_{k,j}^2).\]
Let $({\rm H}_2)$ denote the conditions that the global mean is null, the global variance finite, and for a $p>4$, the centered $p$th moment of the $Y_{k,j}$'s are finite:
\[({\rm H}_2):=\left(
\begin{array}{l}
{\bf m}=0 \text{ and }\beta\in(0,+\infty),\\\hbox{there exist }p>4  \hbox{ s.t. for any } (k,j), 1\leq j \leq k\leq K,  `E\l(| Y_{k,j}-m_{k,j}|^{p}\r)<+\infty.
\end{array}
\right).\]

\begin{figure}[htbp]
\psfrag{v}{$\varnothing$}\psfrag{1}{1}\psfrag{2}{2}\psfrag{3}{3}\psfrag{31}{31}\psfrag{311}{311}\psfrag{11}{11}\psfrag{12}{12}
\centerline{\includegraphics[height=2.5cm]{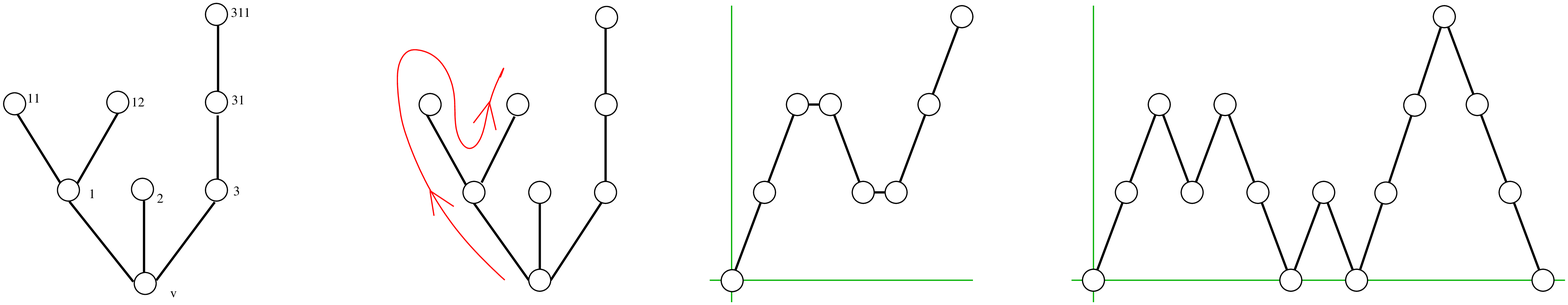}}
\captionn{A tree on which is indicated the depth first traversal, its height and  contour process.}
\end{figure}\par
\subsection*{Encoding of branching random walks}
We study the asymptotic behavior of branching random walks via their encoding by depth first traversal. The depth-first traversal of a tree $ T\in{\cal T}_n$ is a 
function:
\[F_{T}:\{0,..., 2n\}\rightarrow \{\textrm{ vertices of }T\ \},\]
which we regard as a walk around $T$, as follows: 
${F}_{T}(0)=\varnothing$, and  
given ${F}_{T}(i)=z$, choose if possible and according to the LO, the smallest child $w$ of $z$ which has not already been visited, and set $F_T(i+1)=w$. If not
possible, let $F_T(i+1)$ be the father of $z$.\par
We also denote by $\precc$ the lexicographical order (LO) on the planar trees (and $u\prec v$ if $u\precc v$ and $u\neq v$), and let $u(k)$ be the $k$-th vertex
in the LO ($u(0)=\varnothing$). \medskip

We now encode the branching random walk with the help of a pair of processes.
For any $k\in\cro{0,|T-1|}$, let $H^T_k=|u(k)|$ and $R^{T}_k=\ell(u(k))$. The  {height process} $(H^{T}_{s},s\in[0,|T-1|])$ and head label process $(R^{T}_{s},s\in[0,|T-1|])$ are obtained from the sequences  $(H^T_k)$ and $(R^{T}_k)$ by linear interpolation. Alternatively, one may encode the branching random walk with a pair of processes associated with the depth first traversal: for any $k\in\cro{0,2|T|-1}$, let 
$\wh{H}^T(k)=|F_T(k)|$ and $\wh{R}^{T}_k= \ell(F_{T}(k))$. The processes $(\wh{H}^{T}_s,s\in[0,2|T|-1])$ and   $(\wh{R}^{T}_s,s\in[0,2|T|-1])$, obtained by interpolation, are called respectively the contour process and the contour label process; the pair $(\wh{H}^T,\wh{R}^{T})$ is called the head of the discrete snake.
\begin{figure}[htbp]
\psfrag{a}{4}\psfrag{b}{5}\psfrag{5}{5}\psfrag{7}{7}\psfrag{6}{6}\psfrag{9}{9}\psfrag{3}{3}\psfrag{-3}{$-3$}\psfrag{18}{18}\psfrag{R_n}{$R_{9}$}\psfrag{H_n}{$H_{9}$}\psfrag{R_n'}{$\widehat{R_{9}}$}\psfrag{H_n'}{$\widehat{H_{9}}$}\psfrag{1}{1}\psfrag{0}{0}\psfrag{-1}{$-1$}\psfrag{2}{2}\psfrag{4}{4}
\centerline{\includegraphics[height=5.5cm,width =13cm]{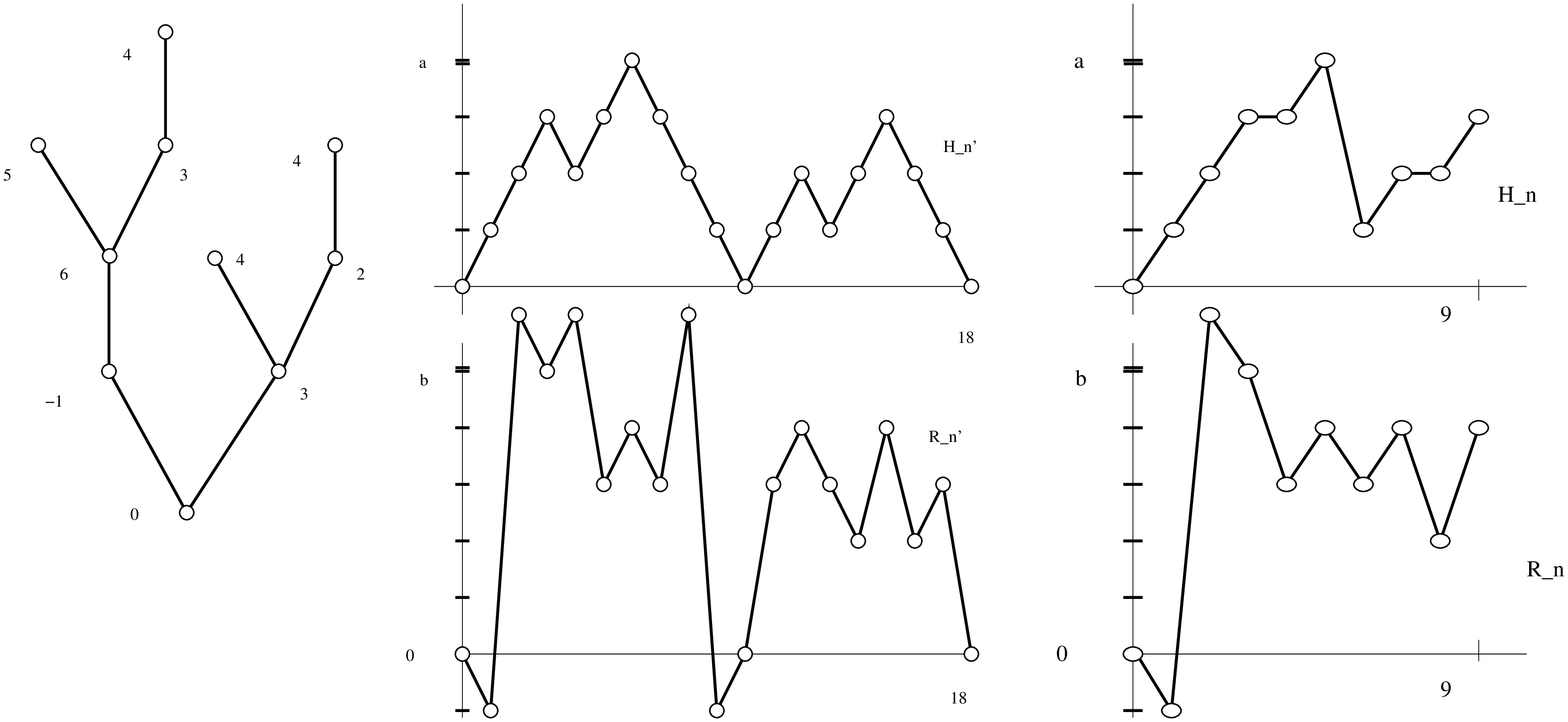}}
\captionn{\label{exm}A branching random walk from ${\cal B}_{9}$. On the first column, the contour process and the  contour label process, on the second column, the height process and the height label process.}
\end{figure}

Let  $\bh_n$, $\widehat\bh_n$,  $\br_n$ and $\widehat{\br}_n$ be the normalized versions of $H^{\bt}$, $\widehat{H}^{\bt}, R^{\bt}$, and $\widehat{R}^{\bt}$ when $\bt$ is  $`P_n$-distributed~:
\[\bh_n(s)=\frac{H_{ns}^\bt}{n^{1/2}},~~~\widehat\bh_n(s)=\frac{\widehat{H}_{2ns}^\bt}{n^{1/2}},~~~\br_n(s)=\frac{R_{ns}^\bt}{n^{1/4}},~~~\widehat\br_n(s)=\frac{\widehat{R}_{2ns}^\bt}{n^{1/4}}, \textrm{ for any }s\in[0,1].\]
Let $d:=\gcd\{k, k\geq 1,\mu_k>0\}$. The support of the distribution of $|\bt|$ -- we write $\supp(|\bt|)$ -- is included in $1+d\,`N$ (and  $`P(|\bt|=1+kd)>0$ for every $k$ large enough).  For $n+1\in \supp(|\bt|)$, the distribution $`P$ under the conditioning by $|\bt|=n+1$ is denoted by $`P_n$, in other words  $`P_n=`P(~.~|\,|\bt|=n+1)$. Even if not recalled, each statement concerning weak convergence under $`P_n$ is assumed to be along the subsequence $(n_k)_k$ for which $`P_{n_k}$ is well defined. In the proofs, we will treat only the case $d=1$, the general case being treated with slight modifications.
\begin{theo}\label{serpdet}
If $({\rm H}_1)$ and $({\rm H}_2)$  are satisfied then  
$$\left(\bh_n,\widehat\bh_n,\br_n, \widehat\br_n\right)\dd
\left(\bh,\bh,\beta{\br},\beta{\br}\right)$$
 in $C([0,1],`R^4)$ endowed with the topology of uniform convergence, where $\bh=2\se/\sigma_{\mu}$ and $\se$ is the normalized Brownian excursion, 
and where conditionally on $\bh$, $\br$ is a centered Gaussian process with covariance function 
\[\cov(\br(s),\br(t))=\check\bh(s,t):=\min_{u\in[s\wedge t, s\vee t]} \bh(u),\textrm{ for any }s,t\in[0,1].\]
\end{theo}

Notice that the same processes $\bh$ and $\br$ appear twice in the limit process. The convergence of processes associated with the contour processes (with a $\widehat{~}$\,) to the same limit as the one associated with the height processes is well understood now, and ``almost'' generic (Duquesne \& Le Gall \cite[Section 2.5]{DUQU1} and \cite{MM3}). In Section \ref{note}, we prove that we may concentrate only of the height process, as done in this paper.
The process $(\br,\bh)$ (or with a different scaling) is called in the literature \it head of the Brownian snake with lifetime process the normalized Brownian excursion (BSBE)\rm . We refer to the works of Le Gall (e.g. \cite{LEG12} and with Duquesne \cite{DG}) for information on the Brownian snake. \par
In this paper, we deal only with the head of the snake, and not precisely in term of snakes, even if, thanks to the homeomorphism theorem \cite{MM3}, evoked below, Theorem \ref{serpdet} has some applications in term of snakes. We refer to \cite{MM3,JM} for the notion of discrete snake which is the discrete analogue of BSBE~: the discrete snake associated with the branching random walk $(T,\ell)$, is the pair $(\widehat{H}^T,\Phi)$ where $\Phi=(\Phi_k)_{k\in\cro{0,2|T-1|}}$ and $\Phi_k$ is the sequence of labels on the branch $\cro{\varnothing,F_T(k)}$. The title of the present paper is then taken from our model of snake under $({\rm H_2})$ in which the global mean is 0.

\subsubsection*{Related works}
The convergence $\widehat\bh_n\dd
\bh$ is due to Aldous \cite{ALD,ALD3} (see also Marckert \& Mokkadem \cite{MM2} for a revisited proof, Pitman \cite[Chap. 5 and 6]{PIT}, and Duquesne \cite{DUQU1} and Duquesne \& Le Gall \cite[section 2.5]{DG} for generalization to GW trees with offspring distribution having infinite variance).\medskip 

The two first results concerning the convergence of discrete snakes to the BSBE appeared in two independent works~:\\
$\bullet$ Chassaing \& Schaeffer \cite{CS} deal with discrete snakes built on underlying trees chosen uniformly in ${\cal T}_n$ (this corresponds to the case where $\mu\sim Geom(1/2)$) and where the displacements are i.i.d., and for any $k,j$, $\nu_{k,j}$ is the uniform distribution in $\{-1,0,+1\}$. They show the convergence of the head of the snake for the Skohorod topology, and the convergence of the moments of the maximum of $\bf r_n$ are also given. This study was motivated by the deep relation between this model of discrete snake and  random rooted quadrangulations, underlined by the authors.\\
$\bullet$ Marckert \& Mokkadem \cite{MM3} studied also the case $\mu\sim Geom(1/2)$ but with  more general centered displacements that have moments of order $6+`e$ (the distribution $\nu_{k,j}$ does not depend on $k,j$, but $\nu_k$ is not assumed to be $\nu_{k,1}\times\dots\times \nu_{k,k}$).  The convergence of the head of the snake holds in $(C[0,1],`R^2)$ and the convergence of the snake itself is given thanks to a ``homeomorphism theorem'' which implies that the convergence of the snake and of its tour (in space of continuous functions) are equivalent. Here it implies that under the hypothesis of  Theorem \ref{serpdet}, the discrete snake associated with our model of labeled trees  converges weakly to the BSBE (see \cite{MM3} for more details).\par 
Then some generalizations appears few months later:\\
$\bullet$ Gittenberger \cite{G1} provides a generalization of a lemma from \cite{MM3} and consider snakes with underlying trees GW trees conditioned by the size (condition equivalent to ${\rm H}_1$). The displacements must be centered and have moments of order $8+`e$. \\
$\bullet$ Janson \& Marckert \cite{JM} show that in the i.i.d. case ($\nu_{k,j}$ do not depend on $(k,j)$),  moments of order $4+`e$ are necessary and needed to get the convergence the BSBE. If no such moment exists the convergence to a ``hairy snake'' is proved under the Hausdorff topology. \\
$\bullet$ In Marckert \& Miermont \cite{GM}, the case of $\nu_{k,j}$ depending of $k,j$ is investigated (also the underlying GW trees are allowed to have two types). The hypothesis are for each $k,j$, $m_{k,j}=0$, condition $({\rm H}_2)$ is satisfied, and then $\sum_{k,j} \mu_k \sigma_{k,j}^2<+\infty$. A motivation was to generalize the works of Chassaing \& Schaeffer \cite{CS} concerning quadrangulations to bipartite maps.\medskip 

Another important point is the convergence of the occupation measure of the head of the discrete snake to the one of the BSBE, the random measure named ISE (the integrated superBrownian excursion introduced by Aldous \cite{ALDMA}, see also Le Gall \cite{LEG12} and \cite{MM3,JM}). Using the convergence of discrete snake to the BSBE, Bousquet-Mélou \cite{B1} and Bousquet-Mélou \& Janson \cite{BJ} deduce new results on ISE and on the BSBE;  for example, some properties on the support of ISE, and of the random density of ISE are derived.  We refer also to Le Gall \cite{LEG4} for the convergence of discrete snake conditioned to stay positive. 

The novelty in the present paper is that the condition $\{m_{k,j}=0, \forall k,j\}$ is replaced by ${\bf m}=\sum_{k\geq 1} \sum_{j=1}^k\mu_k m_{k,j}=0$. This allows to consider some natural models where, for example, the displacements are not random knowing the underlying tree (see Section \ref{exam}). The proof of Theorem \ref{serpdet} relies in part on some results from \cite{GM}, and on a new approach, necessary to control the contribution of the mean of the displacements;  the main point for this, is the comparison of the lineage of each node, with some multinomial r.v. ~: this is the aim of Theorem \ref{zozo}, that we think interesting in itself, since it reveals a thin global behavior of GW trees conditioned by the size.
Unfortunately, the price of this generalization is to consider only offspring distribution with bounded support. The reason comes from the proof of Theorem \ref{zozo}. We guess that some generalization for all families of GW trees (with finite variance) may be found, but for this, a control of an infinite sequence of processes arising in Theorem \ref{zozo} should be provided, what we were unable to do.

\subsection{On the lineage of nodes}

Assume that $({\rm H}_1)$ and $({\rm H}_2)$ holds.
Let $K$ be a bound of the offspring distribution. For $u=i_1\ldots i_h\in \bt$, let 
$u_j=i_1\ldots i_j$ and $\cro{\varnothing,u}=\{\varnothing=u_0, u_1,\ldots, u_{|u|}\}$ be the ancestral line of $u$ back to the root.
Conditionally on $\bt$, $l(u)$ owns the following representations~:
\begin{equation}\label{uuu}\ell(u)=\sum_{m=1}^{|u|} \ell(u_m)-\ell(u_{m-1}).\end{equation}
where $l(u_m)-l(u_{m-1})$ is $\nu_{k,j}$-distributed when $c_{u_{m-1}}(T)=k$ and $i_m=j$, where $\nu_{k,j}$ is the $j$th marginal of $\nu_k$,  and where the r.v.  $(\ell(u_m)-\ell(u_{m-1}))$'s are independent; the variables $\ell(u_m)-\ell(u_{m-1})$ will be often called displacements.

\begin{figure}[htbp]
\psfrag{u}{$u$}\psfrag{hc}{$\wh{H}^{\bf f}$}
\centerline{\includegraphics[height=3cm,width=4cm]{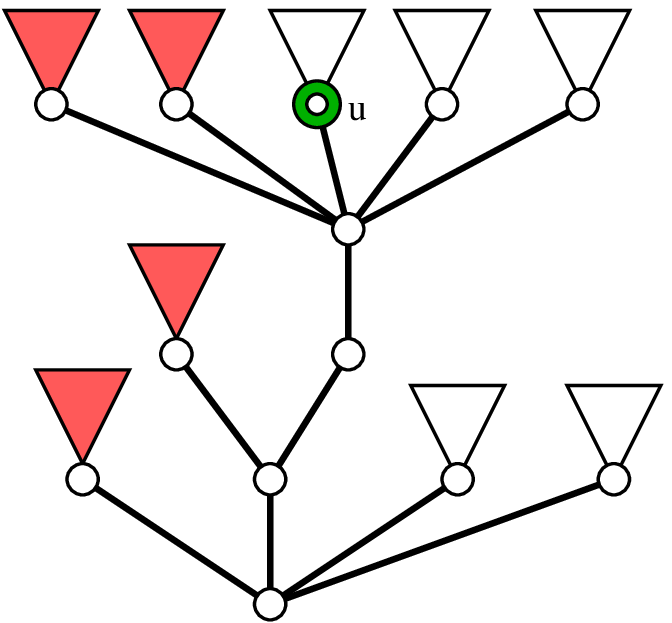}}
\captionn{\label{tret}On this tree $A_{u,1,1}=1, A_{u,2,2}=1,A_{u,4,2}=1, A_{u,5,3}=1$, the others $A_{u,i}$ are 0.} 
\end{figure}

Consider the array $I_K=\l\{(k,j), 1\leq j\leq k\leq K\r\}$. 
Let $u$ be a node of $T$.  For any $(k,j)\in I_K$, let $A_{u,k,j}(T)$ be the number of strict ancestors $v$ of $u$ (the nodes $v \in \llbracket \varnothing,u\llbracket$) such that $c_v({T})=k$, and such that $u$ is  a descendant of $vj$, the $j$th child of $v$ (we write $f_v(u)=j$). We say that $v$ is an ancestor of type $k,j$ of $u$, and we call the vector $A_{u}=(A_{u,i})_{i\in I_K}$ the lineage of $u$ (or the content of $\cro{\varnothing,u}$). See Figure \ref{tret}.
\par
By \eref{uuu}, conditionally on $\bt$, the label  $\ell(u)$ owns the following representations~:
\[\ell(u)\sur{=}{(d)}\sum_{(k,j)\in I_K}\sum_{l=1}^{A_{u,k,j}} Y_{k,j}^{(l)},\]
where the  r.v.   $Y_{k,j}^{(l)}$ are independent, and where for any $l$, $Y_{k,j}$ is $\nu_{k,j}$ distributed. In order to make more apparent the contribution of the $m_{k,j}$'s , and using that ${\bf m}=0$, write
\begin{equation}\label{ellu}
\ell(u)\sur{=}{(d)}\sum_{(k,j)\in I_K}\sum_{l=1}^{A_{u,k,j}} \b(Y_{k,j}^{(l)}-m_{k,j}\b)+\sum_{(k,j)\in I_K}\l(A_{u,k,j}-\mu_k|u|\r)m_{k,j}.
\end{equation}
Assume that $\bt$ is $`P_n$ distributed, and that $u=u(ns)$ for some $s\in(0,1)$. Conditionally on $|u|$, we will see that both parts of the right hand side of \eref{ellu} divided by $n^{1/4}$ converge in distribution, and the limit r.v.  are independent~: in the first part, the fluctuations of $A_{u,k,j}$ around $\mu_k|u|$ are not important, while in the second sum only the fluctuations of $A_{u,k,j}$ around $\mu_k|u|$ matter.

We now concentrate on the r.v.  $(A_{u})'s$ under $`P_n$. 
For any $l\in\cro{0,n}$, $(k,j)\in I_K$, set
\[\bg^{(n)}_{(k,j)}(l):=A_{u(l),k,j}-\mu_k |u(l)|.\] 
For every ${(k,j)\in I_K}$, the process $l\to\bg^{(n)}_{(k,j)}(l)$ encodes the evolution of the number of ancestors of type $k,j$ of $u(l)$, when $l$ varies. 
Consider $\bGn=(\bGn(s))_{s\in[0,1]}$  the process taking its values in $`R^{I_K}$ defined by~: For any $s$, $\bG^{(n)}(s)=\big(\bGn_{k,j}(s)\big)_{(k,j)\in I_K}$ where $s\to\bGn_{k,j}(s)$ is the real continuous process  that interpolates $\bg^{(n)}_{k,j}$ as follows~:
\begin{equation}\label{gkj}
\bGn_{k,j}(s):=\frac{\bg^{(n)}_{k,j}(\floor{ns})+\{ns\}\big(\bg^{(n)}_{k,j}(\floor{ns+1})- \bg^{(n)}_{k,j}(\floor{ns})\big)}{n^{1/4}},~~~s\in[0,1].\end{equation}
The random process $\bGn$ encodes the lineage of all the nodes of $\bt$, and its limiting behavior is described by the following theorem. \par

\begin{theo}\label{zozo}Under $({\rm H}_1)$, $({\rm H}_2)$ the following convergence in distribution holds in $C([0,1])^{\#I_K}\times C[0,1]$ endowed with the topology of the uniform convergence
\[(\bG^{(n)},\bh_n)\dd (\bG,\bh)\]
where  $\bh$ is defined as in Theorem \ref{serpdet} and $\bG=(\bG_{k,j}(s))_{(k,j)\in I_K, s\in[0,1]}$ is a real centered Gaussian field with the following covariance function~:  for any $(k,j)$ and $(k',j')$ in $I_K$, $s$ and $s'$ in $[0,1]$,
\begin{equation}\label{covc}
\cov\big(\bG_{k,j}(s),\bG_{k',j'}(s')\big)= \l(-\mu_k\mu_{k'}+\mu_k\1_{(k,j)=(k',j')}\r)\check{\bh}(s,s').
\end{equation}
\end{theo} 

\subsection{Comments, examples and applications}\label{exam}

1) Theorem \ref{zozo} may be considered as the strongest result of this paper. It gives very precise information on the asymptotic behavior of the process $\bG_n$ that encodes the lineage of all the nodes. This gives a ``global asymptotic'' property reminiscent of the properties of the distinguished branch in ``a size biased GW  tree'' (see \cite[chap. 11]{LY}). The restriction to offspring distributions having a bounded support comes from the proof of this result.

2) For any fixed $(k,j)\in I_K$, knowing $\bh$, $\bG_{k,j}$ is a Gaussian process with  covariance function
\[\cov\big(\bG_{k,j}(s),\bG_{k,j}(s')\big)= \l(-\mu_k^2+\mu_k\r)\check{\bh}(s,s').\]
In other words, the processes $(\bG_{k,j},\bh)$ has the same distribution as $(\sqrt{-\mu_k^2+\mu_k}\,\br,\bh)$, and then up to some multiplicative constants, $(\bG_{k,j},\bh)$ is the head of a BSBE. 
As a simple consequence of Theorem \ref{zozo}, we have 
that $(\bG_{k,j},\bh)_{(k,j)\in I_K}$ is a sequence of heads of BSBE, and that for any  $(k,j)\in I_K$, 
\begin{equation}\b(\bG^{(n)}_{k,j},\bh_{n}\b)\dd\b(\bG_{k,j},\bh\b).\end{equation}The dependence between the different processes $\bG_{k,j}$ is ruled out by \eref{covc}. For any families of real numbers $(\lambda_{k,j})_{(k,j)\in I_K}$, we have
\begin{equation}\l(\sum_{(k,j)\in I_K} \lambda_{k,j}\bG^{(n)}_{k,j},\bh_{n}\r)\dd\l(\sum_{k,j} \lambda_{k,j}\bG_{k,j},\bh\r).\end{equation}

  We would like to stress on the following point: discrete snake are usually constructed with ``two levels of randomness''~: the underlying trees are random and so are the displacements given the underlying tree, and then BSBE appears to be a natural limit of these objects . Here, we provide some objects with only ``one level of randomness'' that converge to the Brownian snake. The BSBE appears as a kind of internal complexity measure in trees measuring the difference between the number of ancestors of type $k,j$ and some expected quantities.

\medskip

3) Consider the case $\mu=\frac{1}{2}(\delta_0+\delta_2)$, $\nu_{2}=\delta_{(+1,-1)}$, of binary trees in which the displacements are not random: $\ell(u1)-\ell(u)=+1$ and $\ell(u2)-\ell(u)=-1$. 
We have ${\bf m}=0$ and $\beta^2=\frac{1}{2} (1+1)=1$ and Theorems \ref{serpdet} and \ref{zozo} apply. Hence, the clear positive bias for $R_n(t)$ for small values of $t$, disappears at the limit.
Note also that this normalizing factor is exactly the same as if  $\nu_{2}=\frac{1}{2}(\delta_{(+1,-1)}+\delta_{(-1,+1)})$ (case where $(\ell(u1)-\ell(u),\ell(u2)-\ell(u))$ is equally likely $(+1,-1)$ or $(-1,+1)$) and as if  $\nu_{2}=(\frac{1}{2}(\delta_{+1}+\delta_{-1}))^2$ (case where the $\ell(u1)-\ell(u)$ and $\ell(u2)-\ell(u)$ are i.i.d., uniform on $\{-1,1\}$).  The question of the convergence of the discrete snake in the case $\nu_{2}=\delta_{(+1,-1)}$ appears first in Marckert \cite{MAR1} in relation with some properties of the rotation correspondence, and the difference between left and right depth in binary trees . The convergence of $(\br_n)$ is not given in \cite{MAR1}, but the convergence of the occupation measure of $\br_n$,  ``the discrete ISE'', to ISE is established. We refer also to Janson \cite{J1} for recent developments concerning the same question.\par
Further, notice that in this model, the label $\ell(u)$ of a vertex $u$ is $\ell(u)=A_{u,2,1}-A_{u,2,2},$
that is the number of left steps minus the number of right steps necessary to climb from the root to $u$ in the binary tree. The convergence of $(\br_n)$ can be seen directly via the one of $(\bG^{(n)})$. 
\begin{equation}\label{GG}
(\bG_{2,1}^{(n)},\bG_{2,2}^{(n)},\bh_n)\dd \b(\bG_{2,1},\bG_{2,2},\bh\b),
\end{equation}
and then $\br_n=\bG_{2,1}^{(n)}-\bG_{2,2}^{(n)}\dd\bG_{2,1}-\bG_{2,2}$ which is, conditionally to $\bh$ and according to \eref{covc}, a centered Gaussian process with covariance function $\check{\bh}(s,t)$. Here, the convergence of $(\br_n)$ appears to be a consequence of the convergence of $\bG_{2,1}$ and $\bG_{2,2}$, encoding the right depth and the left depth in binary trees.  \par

\section{Proofs}
The proofs rely on a precise study of the lineage of the nodes under $`P_n$ and in particular on the comparison of $A_u$ with a multinomial random variable. For this reason we first give some elements on multinomial distributions and on their asymptotic behaviors. We then proceed to the proof of Theorem \ref{zozo}, showing first the convergence of the uni-dimensional distribution then the convergence of the finite-dimensional distribution. The proof of Theorem \ref{serpdet} is given afterward. We think that some points of view especially in the description of the distribution of the lineages in trees under $`P_n$ should provide some new approaches to study the trees under $`P_n$.

\subsection{Prerequisite on multinomial distributions}\label{mere}
The contents of this section is quite classical. 
Consider ${\bf p}=(p_i)_{i\in I_K}$ the distribution on $I_K$, defined by
\[p_{k,j}:=\mu_k ~~~~\textrm{ for any }(k,j)\in I_K.\]
We say that $\m^{(h)}$ is a multinomial r.v. with parameter $h$ and ${\bf p}$, if, for any ${\sf m}=({\sf m}_i)_{i\in I_K}$
\[`Q_h(\{{\sf m}\}):=`P( \m^{(h)}={\sf m})=\binom{h}{({\sf m}_i)_{i\in I_K}}\prod_{i\in I_K} p_i^{{\sf m}_i}\,\,`I_{\ni[h]} ({\sf m})\]
where $\binom{h}{({\sf m}_i)_{i\in I_K}}=h!/(\prod_{i\in I_K}{\sf m}_i!),$ and where for any $n\geq 1$,  $\ni[n]$ is the set of elements $c=(c_{i})_{i\in I_K}$ of ${\mathbb{N}}^{\#I_K}$, such that $ \sum_{i\in I_K}c_{i}=n$. 

Recall that for any $i\in I_K$, ${\cal M}_i^{(h)}$ is a binomial r.v. with parameters $n$ and $p_i$.\par

 In order to fit with further considerations, we introduce the $\#I_K$ dimensional real vector
 ${\cal G}({n,h})=({\cal G}_i({n,h}))_{i\in I_K}$ defined  by
\[{\cal G}_{k,j}(n,h)= n^{-1/4}\big(\m^{(h)}_{k,j}-\mu_{k}\, h\big)~~~\textrm{ for any }(k,j)\in I_K.\]
Let ${\cal G}_{\infty}=({\cal G}_{\infty,_i})_{i\in I_K}$ be a centered Gaussian vector having as covariance function 
\begin{equation}
\cov({\cal G}_{\infty,i},{\cal G}_{\infty,i'})=-p_ip_{i'}+p_i\1_{i=i'} ~~~~\textrm{ for any }i,i'\in I_K.
\end{equation}

\begin{pro}\label{oo}
 Let $(h(n))$ be a sequence of positive integers s.t. $h(n)/\sqrt{n}\to \lambda \in(0,+\infty)$. Under $({\rm H}_1)$ we have  ${\cal G}({n,h(n)})\dd\sqrt{\lambda}\,{\cal G}_{\infty}$ in $`R^{\#I_K}$.
\end{pro}
\proof 
 This may be proved using classical tools. As pointed out by E. Rio in a personal discussion, this is also a consequence of the convergence of the empirical process to the Brownian bridge.
We only sketch the proof (for $\lambda=1$)~:  let $(U_l)_l$ be a sequence of i.i.d. r.v.  uniform on [0,1]. 
Let $F_n$ be the associated empirical distribution function and $F$ the distribution function of $U$. 
Denote by $g_n=F_n-F$. According to Donsker \cite{DO1}, 
$\sqrt{n}g_n\dd {\sf b}$ where ${\sf b}$ is a normalized Brownian bridge. \par
Take ${\bf q}=(q_l)_{l\in\mathbb{N}}$ a distribution on $\mathbb{N}$ and consider ${\cal N}_k^{(n)}=\#\{j, 1\leq j\leq n, U_j\in[q_{1}+\dots+q_{k},q_{1}+\dots+q_{k+1}]\}$. 
Then $({\cal N}_k^{(n)})_{k\geq 1}$ is a multinomial r.v. with parameters $n$ and ${\bf q}$ and satisfies
\[({\cal N}_k^{(n)}-q_k n)/\sqrt{n}=\sqrt{n}\b(g_n(q_1+\dots+q_{k+1})-g_n(q_1+\dots+q_k)\b).\]
By Donsker, for any $L>0$, $\big(({\cal N}_k^{(n)}-q_k n)/\sqrt{n}\big)_{k\leq L}$ converges in distribution to $({\sf b}_{q_1+\dots+q_{k+1}}-{\sf b}_{q_1+\dots+q_{k}})_{k\leq L}$. The properties of ${\sf b}$ allow to conclude.  $\Box$ \medskip

The following Proposition will be used in the proof of the tightness of $(\bG^{(n)})$. 
\begin{pro}\label{qe}
Under $({\rm H}_1)$, for any $\beta>1$,   there exists $c>0$ such that, for any $h>0$, any $n>0$,
\[`E\l(\big\|{\cal G}({n,h})\|_1^\beta \r)\leq c \l({h}/{\sqrt{n}}\,\r)^{\beta/2}.\]
\end{pro}
Recall that all the norms are equivalent in $`R^{\#I_K}$. Here, we use $\|U\|_1=\sum_{(k,j)\in I_K} |U_{k,j}|$.~\\
\proof 
First, since $\|U\|_1^\beta\leq c \sum |U_{k,j}|^\beta$ for some $c>0$,
$`E\l(\big\|{\cal G}({n,h})\|_1^\beta \r)\leq c \sum_{(k,j)\in I_K} `E(|n^{-1/4}\big(\m^{(h)}_{k,j}-\mu_{k}\, h\big)|^{\beta} )$. Since $\m^{(h)}_{k,j}$ is a binomial random variable with parameter $\mu_k$ and $h$, $`E(|\big(\m^{(h)}_{k,j}-\mu_{k}\, h\big)|^{\beta})\leq C({\mu_k,\beta}) h^{\beta/2}$ where the constant $C({\mu_k,\beta})$  depends on $\mu_k$ and $\beta$ (see Petrov \cite{PET}, th. 2.10 p.62).~~$\Box$

\subsection{Decomposition of trees using the lineages}

\begin{figure}[htbp]
\psfrag{u(m)}{$u(m)$}\psfrag{v}{$\varnothing$}
\centerline{\includegraphics[height=3cm, width=2.5cm]{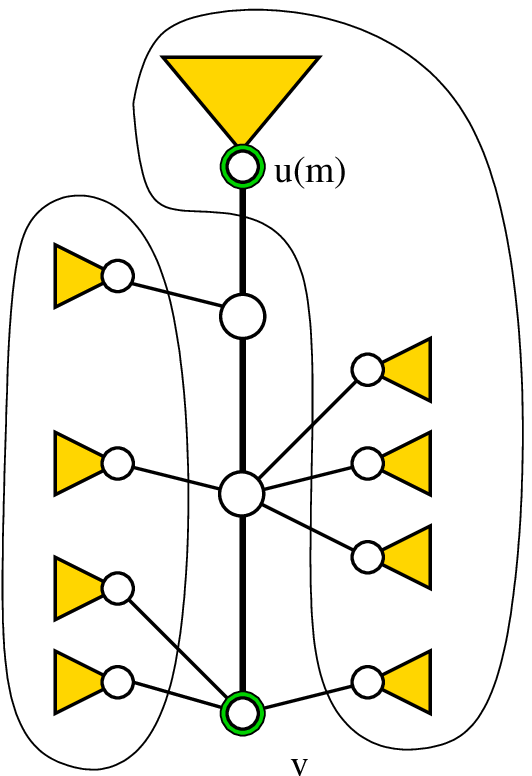}}
\captionn{\label{um} The two forests considered in the decomposition}
\end{figure}

For any $k\in\mathbb{N}$, a forest with $k$ roots is a $k$-tuple of planar trees ${f}=(t^1,\ldots,t^k)$. The size $|f|$
 of ${f}$ is $|t^1|+\ldots|t^k|$.
We denote by ${\bf f}_{k}=(\bt^1,\dots,\bt^k)$ a random forest in which the trees $\bt^1,\dots,\bt^k$ are i.i.d. GW trees with offspring distribution $\mu$. For any $\sa=(\sa_{k,j})_{(k,j)\in I_K}\in {`R}^I$, write
\[N_1(\sa)=\sum_{(k,j)\in I_K} (j-1)\sa_{k,j}~~~~ \textrm{ and }~~~~N_2(\sa)=\sum_{(k,j)\in I_K} (k-j) \sa_{k,j}.\]
\begin{pro}\label{ouch} Let $h$ be a non-negative integer. For any  $\sa\in{\mathbb N}^I[h]$, and any $m\in\cro{0,n}$~:
\begin{equation}
\label{rob}`P_n\l(A_{u(m)}=\sa\r)=`Q_{h}(\sa)\frac{`P\l(|{\bf f}_{N_1(\sa)}|=m-h,\,|{\bf f}\,'_{1+N_2(\sa)}|=n+1-m\r)}{`P(|\bt|=n)},
\end{equation}
where ${\bf f}$ and ${\bf f}\,'$ are two independent forests.
\end{pro}
\proof 
To build a tree $T$ of  ${\cal T}_n$ such that $A_{u(m)}=\sa$, we first build the branch $b=\cro{\varnothing,u(m)}$~: Exactly $\sa_{k,j}$ ancestors $v$ among the $h$ strict ancestors of $u$ satisfy $(c_v(T),f_v(u))=(k,j)$. Hence, there are $\binom n\sa$ way to build  $b$. Then, we complete $b$ in grafting on its neighbors some subtrees satisfying the following constraints. When $A_{u(m)}=\sa$, the number of subtrees rooted on the neighbors of the branch $\llbracket\varnothing,u(m)\llbracket$ visited before $u(m)$ (resp. after $u(m)$) are respectively
\be
N_1(\sa)&=&\#\b\{w, d(\llbracket \varnothing,u(m)\llbracket,w)=1,  w \prec u(m)\b\},\\
1+N_2(\sa)&=&\#\l(\{u(m)\}\cup\b\{w, d(\rrbracket u(m),\varnothing\rrbracket,w)=1, u(m) \prec w\b\}\r).
\ee
See an illustration on Figure \ref{um}.
The $N_1(\sa)$ subtrees must contain exactly $m-|u(m)|$ nodes (the nodes, among the $m+1$ first, not on $\cro{\varnothing, u(m)}$), and the $1+N_2(\sa)$ subtrees must contain exactly $n+1-m$ nodes (the nodes visited after $u(m)$, $u(m)$ included). In other words, we need two forests containing respectively $m-h$ and $n+1-m$ nodes. Hence, using simple considerations on the probability distribution of GW trees we get the announced result. ~~$\Box$ \medskip

A consequence of this Proposition is
\ben\label{decr}
`P_n\l(|u(m)|=h\r)&=&\sum_{\sx \in `N^I[h]}`Q_{h}(\sx) \frac{`P(|{\bf f\,}_{N_1(\sx )}|=m-h,|{\bf f\,}'_{1+N_2(\sx )}|=n+1-m)}{`P(|\bt|=n)}\\
&=&\frac{`P\l(|{\bf f}_{N_1(\m\hh{h})}|=m-h,|{\bf f\,}'_{1+N_2(\m\hh{h})}|=n+1-m\r)}{`P(|\bt|=n)}.
\een
where $\m^{(h)}$ is a multinomial random variable with parameters $h$ and ${\bf p}$.

\subsubsection{Few fact concerning random forests and random trees}
\label{zrf}

 Let $(W_i)_{i\geq 0}$ be a random walk starting from 0 with i.i.d. increments with distribution $(\tilde{\mu}_k)_{k\geq -1}=(\mu_{k+1})_{k\geq -1}$ (that is with increment $\xi-1$, where $\xi$ is $\mu$-distributed). We have
\begin{lem} Assume $({\rm H}_1)$.\\
 (i) \label{fs}(Otter \cite{OTT})  For any $k\geq 1$ and $n\geq k$,
$\dis`P(|\,{\bf f\,}_k|=n)=\frac{k}{n}\,`P(W_n=-k).$\\
(ii) (Central local limit theorem (CLLT) )
\begin{equation}
\sup_{l\in -n+d\mathbb{N}}\l|\frac{\sqrt{n}}{d}\,`P(W_n=l)-\frac{1}{\sqrt{2\pi}\sigma_{\mu}}\exp\l(-\frac{l^2}{2\sigma_{\mu}^2n}\r)\r|\xrightarrow[~~n~~]{} 0.
\end{equation}
(iii) $\sup_{n\geq 0}\sup_{x\geq 0}~x\,`P(W_n=x)\,<+\infty.$
\end{lem} $(i)$ is often called ``conjugation of tree principle'' or ``cyclical lemma'', and  may be found in Pitman \cite[chap 5.1]{PIT} and is usually attributed to Otter, Kemperman or Dvoretzky-Motzkin. \\
$(ii)$ is usually called the central local limit theorem (see Breuillard \cite{BRE} for a state of the art). Recall that $d$ is the span of $\mu$. The support of $W_n$ is included in $-n+d\mathbb{N}=\{u\in\mathbb{Z}, u=-n+di,i\in\mathbb{N}\}$. 
A consequence of $(i)$ and $(ii)$ is that
\begin{equation}\label{ptn}
`P(|\bt|=n)\sim \frac{dn^{-3/2}}{\sqrt{2\pi}\sigma_{\mu}},
\end{equation}
the equivalent being taken along the subsequence where the left hand side is non-null.\\
\bf Proof of $(iii)$: \rm $\sup_{n\geq 0}\sup_{x\geq c\sqrt{n}}\{x\,`P(W_n=x)\}$ is bounded by the Tchebichev inequality. By $(ii)$, $\sup_{x\leq c\sqrt{n}} \sqrt{n}\,`P(W_n=x)\sous{\tend}{n}\frac{d}{\sqrt{2\pi}\sigma_{\mu}}$, then $\sup_{n\geq 0}\sup_{x\leq c\sqrt{n}} \sqrt{n}\,`P(W_n=x)$ is finite.~$\Box$\medskip

The following Lemma controls the maximum increment in the process $H$ under $`P_n$.  
\begin{lem} \label{degdeg} Assume $({\rm H}_1)$. For any $c>0$ then exists $\rho>0$ such that 
\[`P_n\Big(\max_l\l\{\big||u(l+1)|-|u(l)|\big|\r\}\geq \rho \log n\Big)=O(n^{-c}).\]
\end{lem}
\proof  We just sketch the proof that deeply relies  on the conjugation of tree principle. Take  $n+1$ i.i.d. r.v.  $X_1,\dots,X_{n+1}$, $\mu$-distributed. Conditionally on $\sum_{i=1}^{n+1}(X_i-1)=-1$, among  the $n+1$ shifted sequences $(X_{1},\dots,X_{n+1})$, $(X_2,\dots,X_{n+1},X_1),\dots, (X_{n+1},X_1,\dots,X_n)$, exactly one $(X_1^\star,\dots,X_{n+1}^\star)$ corresponds to a sequence $(c_u,u\in T)$ for a tree $T\in{\cal T}_n$ (where the $c_u$ are sorted according the depth first order), and $(X_1^\star,\dots,X_{n+1}^\star)\sur{=}{(d)}(c_u,u\in \bt)$ for $\bt$ under $`P_n$.  \par
The inequality $\big||u(l)|-|u(l+1)|\big|=h>1$ implies that $|u(l+1)|<|u(l)|$, and the deepest common ancestor $v$ of $u(l+1)$ and $u(l)$ has depth $|u(l+1)|-1$. 
Assume that the tree is visited counterclockwise. The nodes in $\croc{v,u(l)}$ are visited consecutively, and each of them has at least one child. Under $`P$, when traversing the tree clockwise (or by symmetry counterclockwise) the gap between two nodes having zero child is a geometrical r.v. $Geom (\mu_0)$ (we work from now on the usual LO order). Denote by $X_1,\dots,X_{n+1}$ i.i.d. random variables $\mu$-distributed and by $G_1,G_2,\dots$ the successive gaps between the zeros.
\[`P\l(\max_{i} G_i \geq \frac{\rho \log n}2\Big|\sum_{i=1}^{n+1}(X_i-1)=-1\r)=O\l(n^{1/2}`P\l(\max_{i\leq n} G_i \geq \frac{\rho \log n}2\r)\r)=o(n^{-c}),\]
for $\rho$ large enough. Note that the first maximum is taken on a random number of terms, a.s. bounded by $n$. By the conjugation  of tree principle, we get the result.~$\Box$
\begin{remi}\label{rzer}Using the same argument, one may control the depth of the last node $u(n)$~: for any $c>0$ then exists $\rho>0$ such that 
\begin{equation}\label{compl}
`P_n(|u(n)|\geq \rho \log n)=O(n^{-c}).
\end{equation}
\end{remi}
\medskip 

For $u\in T, l\in \cro{0,|u|}$ and $(k,j)\in I_K$, let $A_{u,l,k,j}$  be the number of ancestors 
$v\in\llbracket\varnothing,u\llbracket$ such that $d(u,v)\leq l$, and 
for which $c_u(T)=k$ and $f_v(u)=j$. 

\begin{lem}\label{aukj} $(i)$ For every $c>0$, there exists $\gamma>0$, such that for $n$ large enough,
\[`P_n\left(\exists (k,j)\in I_K, u\in T,
\left|A_{u,k,j}-\mu_k|u| \right| \geq \gamma\sqrt{|u|\log n}\right)\leq n^{-c}\]
$(ii)$ \label{mieux}
For every $c>0$, there exists $\gamma>0$ such that, for $n$ large enough
\[`P_{n}\left(\exists(k,j)\in I_K, u\in T,
 l\in(0,|u|],
\left|A_{u,l,k,j}-\mu_kl\right|\geq \gamma\sqrt{l\log n}\right)\leq n^{-c}.\]
\end{lem}
\proof $(ii)$ clearly implies $(i)$. But let us prove $(i)$ first. Using \eref{rob} and \eref{ptn}, we have for some constant $c>0$, for any $m\in\cro{0,n}$, any $h\geq 1$, any $\sa \in\ni[h]$,
\begin{equation}
`P_n(A_{u(m)}=\sa)\leq cn^{3/2}`Q_h(\sa)\1_{h\leq n}.
\end{equation}
Then $`P_n\l(\exists m \in\cro{0,n}, (k,j)\in I_K, \big|A_{u(m),k,j}-\mu_k |u(m)|\big|\geq \gamma\sqrt{|u(m)|\log n}\r) \leq $
\[cn^{3/2}\sum_{m=0}^n\sum_{h=0}^n`P\Big(\exists(k,j)\in I_K, \big|{\cal M}^{(h)}_{k,j}-\mu_k |h|\big|\geq \gamma\sqrt{ h\log n}\Big).\]
This latter probability is smaller, for any $m\leq n$, $h\leq n$ than $\#I_K n^{-\gamma^2/2}$ by Hoeffding. Hence $`P_n\left(\exists (k,j)\in I_K, u\in T, \left|A_{u,k,j}-\mu_k|u| \right| \geq \gamma \sqrt{|u|\log n}\right)\leq cn^{7/2}n^{-2{\gamma}^2}$.\par
\begin{figure}[htbp]
\psfrag{u}{$u(m)$}
\psfrag{v1}{$v_1$}
\psfrag{v2}{$v_2$}
\psfrag{3}{3}\psfrag{31}{31}\psfrag{32}{32}
\centerline{\includegraphics[width=6cm]{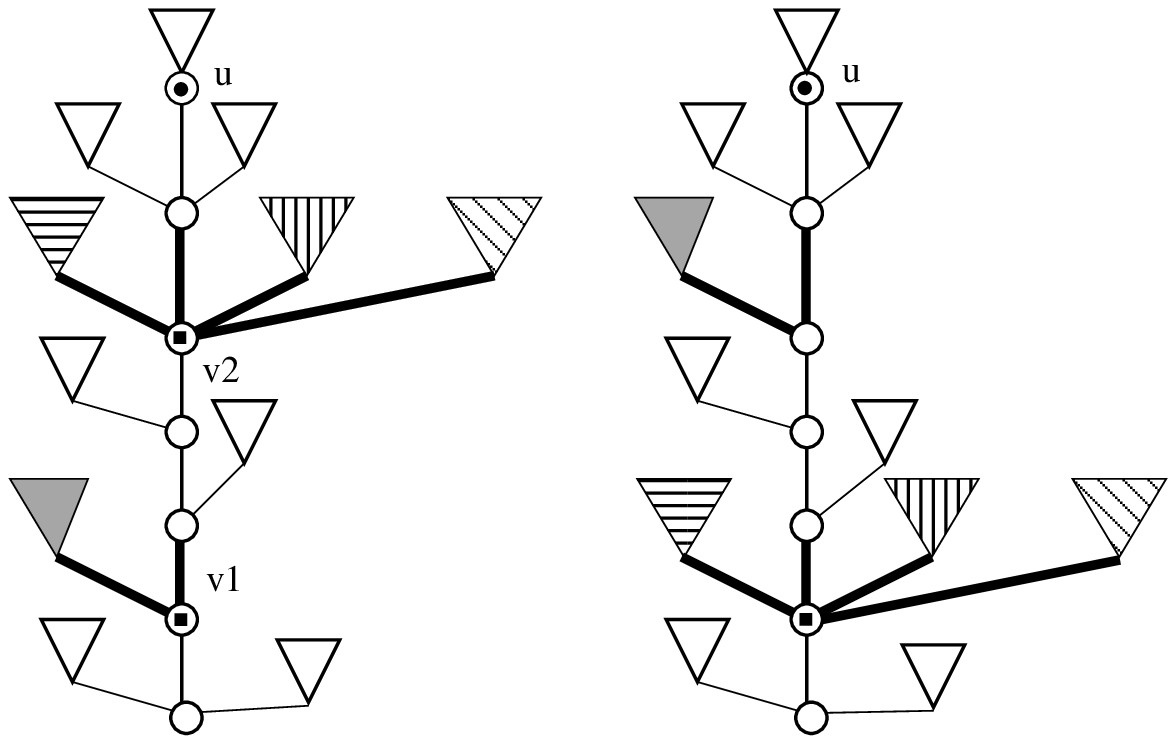}}
\captionn{\label{exch} Exchange of two nodes in a lineage}
\end{figure}

For $(ii)$, assume that $u(m)=h$ and for $l\leq h$, take $v_1,\dots,v_l$ the ancestors of $u(m)$ at depth $0\leq h_1<\dots<h_l<h$, and set $A'_{u(m),l,k,j}=\#\{i, c_{v_i}=k,f_{v_i}(u(m))=j\}$, the lineage of $u(m)$ restricted to the nodes $v_i$'s. By ''symmetry``, $(A'_{u(m),l,k,j})_{k,j}$ and $(A_{u(m),l,k,j})_{k,j}$ have the same distributions. 
Here ``symmetry'' means the following~: let $v_1$ and $v_2$ be two ancestors of $u(m)$. Exchange in $T$, the two nodes $v_1$ and $v_2$ together with the subtrees rooted on their children not on $\cro{0,u(m)}$, as on figure \ref{exch}. We get $T'$. First $T'$ and $T$ has the same weight under $`P_n$. Second, $A_{u(m)}$ has the same value in $T$ and $T'$, and the nodes $u(m)$ in $T$ and $T'$ have the same depth ($u(m)$ is by definition the $m$th node).\par
Now take $v$ the ancestor of $u(m)$ at depth $l$. By symmetry,  $(A_{v,k,j})_{k,j}$ and $(A_{u(m),l,k,j})_{k,j}$ have the same distributions. And thus, by $(i)$, for any $m\leq n$, $l\leq n$, $`P(\exists (k,j)\in I_K,\l|A_{u(m),l,k,j}-\mu_k|l| \right| \geq \gamma\sqrt{l\log n)}$ is certainly smaller than $cn^{7/2}n^{-2{\gamma}^2}$. As a direct consequence, $cn^{7/2+2}n^{-2{\gamma}^2}$ is a bound for $`P_{n}\left(\exists(k,j)\in I_K, u\in T,
 l\in(0,|u|), \left|A_{u,l,k,j}-\mu_kl\right|\geq \gamma\sqrt{l\log n}\right)$.~ ~$\Box$\medskip

We end this section with a result concerning multinomial random variables.
For any $h>0$, set  
\[J_{h}=\l\{\sa \in \ni[h], (N_1(\sa),N_2(\sa))\in \Big[\frac{\sigma^2 h}{2}- h^{2/3},\frac{\sigma^2 h}{2}+ h^{2/3}\Big]^2\r\}.\]
\begin{lem}\label{betat} For  $h\in\mathbb{N}$, $N_1(\m\hh{h})$ and $N_2(\m\hh{h})$ have the same law, and there exists $c_1>0,c_2>0$, s.t
\[`P\l(\m^{h}\notin J_{h}\r)\leq c_1\exp(-c_2\,h^{-1/3}).\]
\end{lem}
\proof The first assertion is easy.
Writing $\{|N_1(\m^{h})-\frac{\sigma^2h}{2}|\geq n^{2/3}\}\subset\bigcup_{k,j}\{|\m^h_{k,j}-h\mu_k|\geq \mu_k h^{2/3}\}$ (recall that $\sum k\mu_k=\frac{\sigma^2}2$), by Hoeffding, one has $`P(|\m^h_{k,j}-h\mu_k|\geq \mu_k h^{2/3})\leq 2\exp(-\mu_k^2/h^{1/3})\1_{\mu_k>0}$. Summing this for $(k,j)\in I_K$, one gets the results.$~\Box$

\subsubsection{A first comparison Lemma}\label{fcr}
 In this section $S$ denote a Polish space. 
For any r.v. $X$ taking its values in $S$, we denote by $`P_X$ the distribution of $X$: that is $`P_X(A)=`P(X\in A)$ for any $A$ Borelian of $S$. 
\begin{defi}\label{def1}
Let $(Y_1,Y_2,\dots)$ and $(X_1,X_2,\dots)$ be two sequences of r.v. taking their values in $S$ such that $`P_{X_n}$ is absolutely continuous with respect to $`P_{Y_n}$, we write $`P_{X_n}\prec `P_{Y_n}$. Let $f_n$ be a negative measurable function $f_n$ such that $`P_{X_n}=f_n `P_{Y_n}$ (which existence is ensured by the Radon-Nikodym theorem)~: for any Borelian $A$ of $S$,  $`P_{X_n}(A)=\int_A f_n d`P_{Y_n}$.  
We say that $`P_{X_n}/`P_{Y_n}\to 1$, or  $X_n\di Y_n\to 1$, if $f_n$ goes to 1 in the following (weak) sense~: for any $`e>0$, the set  $A_{`e}^n=\{x, |f_n(x)-1|<`e\}$ satisfies $`P_{Y_n}(A_{`e}^n)\to 1$.
\end{defi}
If ${X_n}\di Y_n\to 1$ then  $`P_{X_n}(A_{`e}^n)\to 1$, and for any $B\subset  A_{`e}^n$, $|`P(Y_n\in B)-`P(X_n\in B)|\leq `e`P(Y_n\in B)$ and then the total variation distance between $X_n$ and $Y_n$, defined by $\sup_{B \textrm{ Borelian}}|`P(X_n\in B)-`P(Y_n\in B)|$ goes to 0. Hence, the following Lemma is a straightforward  consequence of the Portmanteau theorem~:
\begin{lem}\label{dr}If ${X_n}\di Y_n\to 1$ and $Y_n\dd Y$ then $X_n\dd Y$.
\end{lem}

\subsubsection{Proof of the convergence of the uni-dimensional distributions in Theorem \ref{zozo}}
In this section we work under $`P_n$. 
Let $X_m^n:=(A_{u(m)},|u(m)|)$ and $Y_m^n:=(A^\star_{m},|u(m)|)$ where  the distribution of $A^\star_{m}$ knowing $|u(m)|=h$ is simply $`Q_{h}$.
The aim of this section is to compare $X_m^n$ with $Y_m^n$ and to deduce from the asymptotic behavior of  $Y_m^n$ some information on $X_m^n$. The proof of the convergence of the finite-dimensional distributions will  also use this strategy. \par
For $M>0$, and $n\in\mathbb{N}$, consider
\[\Lambda_{n,M}=\l\{(\sa,h), h\in\sqrt{n}[M^{-1},M], \sa\in J_{h}\r\}.\]
We have
\begin{pro}\label{RN}
$i)$ For any $m,n$, $`P_{X_m^n}\prec `P_{Y_m^n}$.\\
$ii)$ For any $s\in(0,1)$, $\alpha>0$, there exists $M$ s.t. for $n$ large enough,
$`P_n\big(Y_{\floor{ns}}^n\in \Lambda_{n,M}\big)\geq 1-\alpha$\\
and for any $M>0$,
\begin{equation}\label{eqae}
\sup_{(a,h)\in \Lambda_{n,M}}
\l|\frac{`P_n(X_{\floor{ns}}^n=(\sa,h))}{`P_n(Y_{\floor{ns}}^n=(\sa,h))}-1\r|\xrightarrow[n]{~~~}0
\end{equation}
$iii)$ For any $s\in(0,1)$, ${{X_{\floor{ns}}^n}}\di {Y_{\floor{ns}}^n}\to 1.$
\end{pro}

\proof $(iii)$ is a consequence of $(ii)$. 
  Let $\sa\in \mathbb{N}_I[h]$. Since $\{A_{u(m)}=\sa\}\subset \{|u(m)|=h\}$,  $`P_n((A_{u(m)},|u(m)|)=(\sa,h))=`P_n(A_{u(m)}=\sa).$ According to Proposition \ref{ouch}, and Formula \eref{decr}
\ben \label{pass}
\frac{`P_n(X_m^n=(\sa,h))}{`P_n(Y_m^n=(\sa,h))}&=&\frac{`P(|{\bf f}_{N_1(\sa)}|=m-h,\,|{\bf f}'_{1+N_2(\sa)}|=n+1-m)}{`P\l(|{\bf f}_{N_1(\m^{(h)})}|=m-h,|{\bf f}'_{1+N_2(\m^{(h)})}|=n+1-m\r)}.
\een
Then $(i)$ holds true. 
Assume now that $s\in(0,1)$ and $\alpha>0$ are fixed. There exists $M$ such that for $n$ large enough, $`P_n({|u(\floor{ns})|}\in{\sqrt{n}}[M^{-1},M])\geq 1-\alpha/2$ (since $\bh_n\dd\frac{2}{\sigma}{\bf e}$ and since $`P({\bf e}_s=0)=0$ for any $s\in(0,1)$).  
For such a $M$, 
\be
`P_n\l(Y_{\floor{ns}}^n\in \Lambda_{n,M}\r)&=&`P_n\l(Y_{\floor{ns}}^n\in \Lambda_{n,M},|u(\floor{ns})|\in{\sqrt{n}}[M^{-1},M]\r)\\
&=& \sum_{l\in\sqrt{n}[M^{-1},M]} `P(|u(\floor{ns})|=l)`P_n\l(Y_{\floor{ns}}^n\in \Lambda_{n}\,\b|\, |u(\floor{ns})|=l\r)\\
&\geq & `P_n\l(|u(\floor{ns})|\in{\sqrt{n}}[M^{-1},M]\r)\min_{l\in\sqrt{n}[M^{-1},M]}  `P\l(\m^{(l)}\in J_{l}\r) 
\ee
This infimum goes to 1 thanks to Lemma \ref{betat}. 
According to Lemma \ref{fs} $(i)$ and $(ii)$, since $\bf f$ and $\bf f'$ are independent, $`P(|{\bf f}_{N_1(\sa)}|=\floor{ns}-h,\,|{\bf f}'_{1+N_2(\sa)}|=n+1-\floor{ns})=$
\[\frac{N_1(\sa)(1+N_2(\sa))}{(\floor{ns}-h)(n+1-\floor{ns})}`P(W_{\floor{ns}-h}=-N_1(\sa))`P(W_{n-\floor{ns}+1}=-N_2(\sa)-1)\]
and then for any $M>0$, 
\[\sup_{(a,h)\in \Lambda_{n,M}}
\l|\frac{`P(|{\bf f}_{N_1(\sa)}|=\floor{ns}-h,\,|{\bf f}'_{1+N_2(\sa)}|=n+1-\floor{ns})}{q_{n,s,h}}-1\r|\xrightarrow[n]{~~~}0 \] 
for \[q_{n,s,h}=\frac{\sigma^2 h^2\exp\l(-\frac{\sigma^4 h^2}{8{ns}(1-s)}\r)}{8\pi n^3(s(1-s))^{3/2}}.\]
Now,  $`P\l(|{\bf f}_{N_1(\m\hh{h})}|=\floor{ns}-h,|{\bf f}'_{1+N_2(\m\hh{h})}|=1+n-\floor{ns}\r)=A_h+B_h$ where
\be
A_h&:=&`P\l(|{\bf f}_{N_1(\m^{(h)})}|=\floor{ns}-h,|{\bf f}'_{1+N_2(\m^{(h)})}|=1+n-\floor{ns},\m^{(h)}\notin J_h\r)\\
B_h&:=&`P\l(|{\bf f}_{N_1(\m^{(h)})}|=\floor{ns}-h,|{\bf f}'_{1+N_2(\m^{(h)})}|=1+n-\floor{ns},\m^{(h)}\in J_h\r)
\ee
Using again Lemma \ref{fs} $(i)$ and $(ii)$, we get
\[\sup_{h\in\sqrt{n}[M^{-1},M]} \l|\frac{B_h}{q_{n,s,h}}-1\r|\xrightarrow[n]{~~~} 0.\]On the other hand, $A_h\leq `P(\m^{(h)}\notin J_h)\leq c_1\exp(-c_2 h^{1/3})\leq 2\exp(-c\,n^{1/6}/M)$ for any $h\in\sqrt{n}[M^{-1},M]$.
To complete the proof of $(ii)$, check that $\sup_{h\in\sqrt{n}[M^{-1},M]}\l|{A_h}/{B_h}\r|\xrightarrow[n]{~~~} 0.~\Box$

\begin{cor}
\label{opieds} For any $s \in (0,1)$, let $s_n={\floor{ns}}/n$, we have
\[\l(\bGn(s)(s_n),\bh_{n}(s_n)\r)\di\l({\cal G}\l(n,\sqrt{n}\,\bh_{n}(s_n)\r),\bh_n(s_n)\r)\to 1,\]
and the convergence of the uni-dimensional distributions holds in Theorem \ref{zozo}.
\end{cor}
Recall that ${\cal G}$ is defined in Section \ref{mere}.\\
\proof  Proposition \ref{RN} yields the first assertion of the Corollary. For the second one, by Lemma \ref{dr}, it suffices to establish that for any $s\in[0,1]$
\begin{equation}\label{mar}
\l({\cal G}\l(n,\sqrt{n}\,\bh_{n}(s_n)\r),\bh_n(s_n)\r)\dd ({\cal G}_{\infty}^{\bh_s},\bh_s)\sur{=}{(d)}(\bG(s), \bh_s).
\end{equation}
For $s=0$ or $s=1$, this is a consequence of $\bh_n(s)\proba 0$.
For $s\in(0,1)$, since $\bh_n\dd \bh$ in $C[0,1]$, by the Skohorod representation theorem \cite[Theorem 3.30]{KAL}, there exists a probability space  on which this convergence is a.s.. 
On this space (or on an augmented space on which the pair $\l({\cal G}\l(n,\sqrt{n}\,\bh_{n}(s_n)\r),\bh_n(s_n)\r)$ is defined), \eref{mar} holds a.s.. To prove that the convergence of the uni-dimensional distribution  holds in Theorem \ref{zozo}, it remains to control the distance between $\l(\bGn(s_n),\bh_{n}(s_n)\r)$ and $\l(\bGn(s),\bh_{n}(s)\r)$.
Since $\bh_{n}\dd \bh$, $|\bh_{n}(s_n)-\bh_{n}(s)|\proba 0$. For $\bGn$ this is more complex, and we will establish some bounds useful also for the tightness. Let 
\[\Omega_n^{\rho}=\Big\{T\in{\cal T}_n, \max_{l} \b||u(l+1)|-|u(l)|\b|\leq \rho \log n \Big\}.\]
Let $`e>0$. 
According to Lemma \ref{degdeg}, for $\rho$ large enough, $`P_n(\Omega_n^{\rho})>1-`e$ for $n$ large enough.
We have for $s'_n={\floor{ns+1}}/{n}$, 
\ben \label{int}
\|\bGn(s_n)-\bGn(s)\|_1`I_{\Omega_n^{\rho}}&=& n(s-s_n)\sum_{i\in I_K}`I_{\Omega_n^{\rho}}\l|\bGn_{i}(s'_n)-\bGn_{i}(s_n)\r|.
\een
In $\Omega_n^{\rho}$,  the differences $|\bGn_{k,j}(s'_n)-\bGn_{k,j}(s_n)|$ are bounded by $\rho n^{-1/4}\log n$. Hence, since $s-s_n\leq 1/n$, for any $`e'>0$, for $n$ large enough 
\begin{equation}\label{mau}
\|\bGn(s_n)-\bGn(s)\|_1`I_{\Omega_n^{\rho}}\leq c(s-s_n)^{1/4-`e'}.
\end{equation} for some constant $c$.
One concludes that $\|\bGn(s_n)-\bGn(s)\|_1\proba 0$. $~\Box$

\subsection{Convergence of the finite-dimensional distributions}

In this Section, $\kappa\geq 2$ is a fixed integer.  We denote by $\sk$ the vector $(s_1,\dots,s_\kappa)$ where  $0<s_1<\dots<s_{\kappa}\leq 1$ are fixed. Let $T\in{\cal T}_n$. For $i\in\cro{1,\kappa}$, set $u_i=u(\floor{ns_i})$, $u_0=u_{\kappa+1}=\varnothing$. The aim of this section is to study the joint distribution of $(A_{u_i})_{i\in\cro{1,\kappa}}$ under $`P_n$. The ideas are of the same type as in the case of the uni-dimensional distributions, but the details are more involved since the dependences between the r.v. $A_{u_i}$'s must be taken into account. For this, we must consider the shape of the tree spanned by the $u_i$'s.\par
Denote by $\check u_{i,j}$ the deepest (i.e. youngest) common ancestor between $u_i$ and $u_j$.
Let $T_{\sk}=\bigcup_{i=1}^{\kappa} \cro{\varnothing, u_i}$ be the subtree ``spanned'' by the $u_i$'s,
 $L(T)$ the set $\{u_{i},  i \in \cro{0,\kappa}\}$, and $Z(T)=\{\check u_{i,j}, 1\leq i<j\leq \kappa\}=\{\check u_{i,i+1}, i\in\cro{1,\kappa-1}\}$, the set of \it branching nodes \rm in $T_{\sk}$.

\begin{figure}[htbp]
\psfrag{emp}{$\varnothing$}\psfrag{u0}{$u_0$}\psfrag{u1}{$u_1$}\psfrag{u2}{$u_2$}\psfrag{u3}{$u_3$}\psfrag{u4}{$u_4$}\psfrag{u5}{$u_5$}\psfrag{u35}{$\check{u}_{3,5}$}\psfrag{u12}{$\check{u}_{1,2}$}\psfrag{u23}{$\check{u}_{2,3}$}
\centerline{\includegraphics[height=2.7cm]{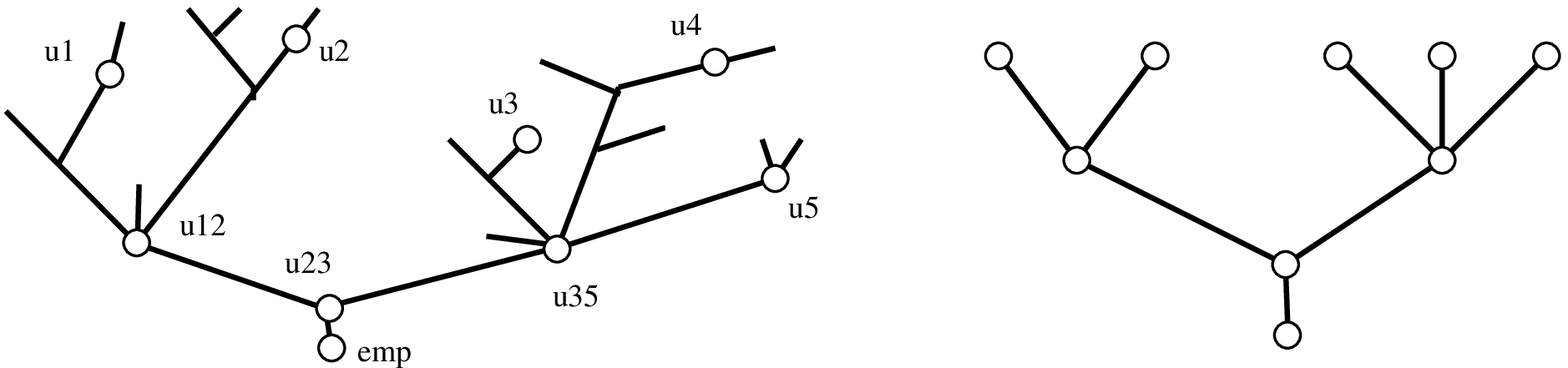}}
\captionn{\label{anc33} A tree $T$ and the associated tree $\Phi(T)=\{\varnothing,1,11,111,112,12,121,122,123\}$.}
\end{figure}

\begin{defi}
The shape function $b$ associates with $T_{\sk}$ the smallest tree having the same shape. Formally
\[\app{b}
{{\cal T}\times [0,1]^{\kappa}}
{{\cal T}}{(T,\sk)}
{T^b}
\]
where $T^b$ is characterized by~:\\
$i)$ $T^b\in{\cal T}$ and $\#T^b=\#\b(Z(T)\cup L(T)\b),$\\
$ii)$ there exists an increasing function $\Phi_{T}$ from $Z(T)\cup L(T)$ in $T^b$, preserving the descendants~: $\Phi_{T}(u)$ is an ancestor of $\Phi_{T}(v)$ in $T^b$ ifff $u$ is an ancestor of $v$ in $T$.
\end{defi}
In other words, $T^b$ is the only tree with $\#\b(Z(T)\cup L(T)\b)$ nodes having the same branching structure as $T_{\sk}$ (see Figure \ref{anc33}); it can be constructed in somehow squeezing the paths between the nodes of $Z(T)\cup L(T)$ in unit length edge (and in renaming the vertices in order to get a tree). The function $\Phi_{T}$ is unique and, for short, for any $u\in Z(T)\cup L(T)$, we write $u^b$ instead of $\Phi_{T}(u)$.

The set $\frak{E}(T)$ of ``spanned branches'' between the nodes of $L(T)\cup Z(T)$ is defined by
\be
\frak{E}(T)&=&\{(u,v), u,v\in L(T)\cup Z(T), u\neq v, u\prec v, \croc{u,v} \cap (L(T)\cup Z(T))=\emptyset\}\\
&=&\{(u,v), u,v\in L(T)\cup Z(T), u^b=\fa(v^b)\},
\ee
where $\fa(u)$ stands for the father of $u$.
Notice that if $(u,v)\in\frak{E}(T)$ then $u$ is an ancestor of $v$. 
The LO order induces an order on  $\frak{E}(T)$ (also denoted by  $\prec$)~: for $e_1=(u_1,v_1)$ and $e_2=(u_2,v_2)\in \frak{E}(T)$,  $e_1\prec e_2$ if $v_1\prec v_2$. We denote by ${\cal H}_T=(\#\croc{u,v})_{(u,v)\in \frak{E}(T)}$
the ordered list of the spanned branches  distances.

In order to control the dependences between the $A_{u_i}$'s, we study analogous quantities associated with spanned branches. For any $(u,v)\in \frak{E}(T)$,  define $A_{(u,v)}$, the content of the edge $(u,v)$ by
\[A_{(u,v),k,j}:=\#\l\{w\in \croc{u,v}, c_w=k, f_w(v)=j\r\}.\]
The contribution of the extremities of the spanned branches are not counted in any of the $A_{(u,v),k,j}$'s in order to simplify the enumerations in the rest of the paper. It is easy to check that 
\begin{equation}\label{eqdiff}
A_{(u,v),k,j}=(A_{v,k,j}-A_{u,k,j})-\1_{(c_u,f_u(v))}(k,j).
\end{equation}
The contributions of the nodes of  $Z(T)$ are encoded by the sequence $\Theta_T$, ordered by $\prec$~:
\[\Theta_{T}=\b(c_u,f_u(L(T))\b)_{u\in Z(T)\cup\{\varnothing\}}.\]
Notice that $f_u(L(T))$ is a subset of $\cro{1,\dots,c_u}$ with $c_{u^b}$ elements.

\subsubsection{Subtrees visited between two elements of $L(T)$}
\label{subs}
\begin{figure}[htbp]
\psfrag{u0}{$u_0$}\psfrag{u1}{$u_1$}\psfrag{u2}{$u_2$}\psfrag{u3}{$u_3$}\psfrag{z1}{$z_1$}\psfrag{z2}{$z_2$}\psfrag{s0}{$F_0$}\psfrag{s1}{$F_1$}\psfrag{s2}{$F_2$}\psfrag{s3}{$F_3$}
\centerline{
\includegraphics[height=6cm]{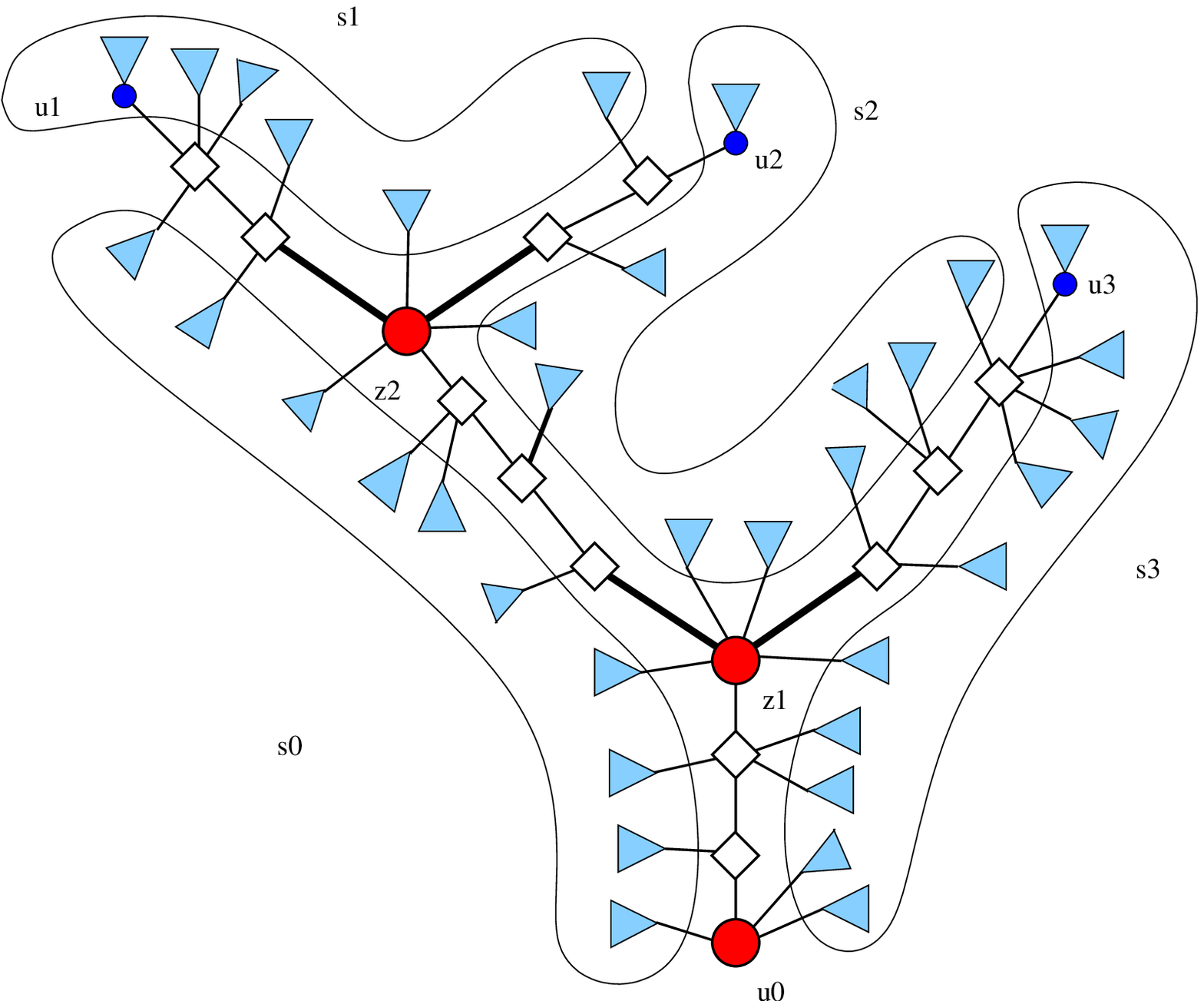}}
\captionn{\label{anc3} The ordered spanned branches are $(\varnothing, z_1), (z_1, z_2), (z_2, u_1), (z_2,u_2), (z_1,u_3)$, $\Theta_T=\big((4,2), (6,2,5),(5,2,4)\big)$. The triangles are subtrees not drawn.  The forests considered in the decomposition are surrendered.}
\end{figure}

(An illustration of the quantities considered in this section is given on Figure \ref{anc3}). 
We denote by
\be
\Sub_0&=&\l\{v\in T, d(v,\llbracket \varnothing, u_1\llbracket)=1, \varnothing \prec v \prec u_1 \r\},\\
\Sub_i&=&\l\{v\in T, d(v,\croc{u_i,u_{i+1}})=1, u_i\prec v \prec u_{i+1}\r\}\cup \{u_i\}~~ \textrm{ for }i\in\cro{1,\kappa},\\
\Sub_{\kappa+1}&=&\l\{v\in T, d(v,\rrbracket u_{\kappa},\varnothing\rrbracket)=1,   u_\kappa \prec v\r\}\cup\{u_{\kappa}\}
\ee
the set of the roots of the subtrees, rooted on the neighbors of $\cro{u_i,u_{i+1}}$, visited by the depth first traversal between $u_i$ and $u_{i+1}$ (up to the borders effects). 
The cardinalities of the set $\Sub_l's$ are characterized by the triplet 
$\left({\cal A}_{T}, \Theta_T,T^b\r),$ where
$ {\cal A}_{T}=\l(A_{(u,v)}\r)_{(u,v)\in \frak{E}_{T}}.$
 For any $l\in \cro{0,\kappa}$, we have
 \begin{equation}\label{tere}
\#\Sub_l(T^b,{\cal A}_T,\Theta_{T})={\cal N}_{l,{l+1}}(T^b,{\cal A}_{T})+{\cal Y}_{l,{l+1}}(T^b,\Theta_T)
 \end{equation}
where, for any $l$,
\[{\cal Y}_{l,{l+1}}(T^b,\Theta_T)=\1_{l\neq 0}+f_{\check u_{l,l+1}}(u_{l+1})-f_{\check u_{l,l+1}}(u_l)-1+\sum_{z^b\in \rrbracket \check u_{l,l+1}^b, u_{l}^b\llbracket} c_{z}-f_z(u_l) +\sum_{z^b\in\rrbracket \check u_{l,l+1}^b,u_{l+1}^b\llbracket}(f_z(u_l)-1)\]
and 
\[{\cal N}_{l,{l+1}}(T^b,{\cal A}_{T})=\sum N_1(A_{(z_1,z_2)})+ 
\sum N_2(A_{(z_1,z_2)}),\]
where the first sum is taken on the pairs $(z_1,z_2)$ such that $(z_1^b,z_2^b)\in\cro{\check u_{l,l+1}^b,u_{l+1}^b}$, and $z_1^b=\fa(z_2^b)$ and the second one, on the pairs $(z_1,z_2)$ such that $(z_1^b,z_2^b)\in\cro{\check u_{l,l+1}^b,u_{l+1}^b}$ and $z_2^b=\fa(z_1^b)$.

This is similar to Proposition \ref{ouch}~: ${\cal N}_{l,{l+1}}(T^b,{\cal A}_{T})$ counts the number of subtrees rooted on the neighbors of the spanned branches , on their right or on their left;  ${\cal Y}_{l,{l+1}}(T^b,\Theta_T)$ counts the number of subtrees rooted on the neighbors of the nodes of $Z(T)$. We let $T_u=\{v\in T_{\infty}:uv\in T\}$ be the {\em fringe subtree} of $T$ rooted at $u$.
The cardinality $F_l=\#(T_u)_{u\in \Sub_l}$ of the forests constituted with the fringe subtrees of $T$ rooted on the nodes of $\Sub_l$ satisfies
\begin{equation}
F_l({\cal H}_T,T^b)= (ns_{l+1}-ns_{l}+1)-(|u_{l+1}|-|\check u_{l,l+1}|)-\1_{l=0},
\end{equation}
since the visit times of the nodes $u_i$ are $ns_i$ and since 
$|u_{l+1}|-|\check{u}_{l,l+1}|+1$ nodes visited during $\cro{ns_l, ns_{l+1}}$ are not in $(T_u)_{u\in \Sub_l}$.

\subsubsection{Decomposition of a tree $T$ given the ${\cal A}_{T}$ and $T^b$} 

Let ${\cal T}^B_{2\kappa-1}=\l\{T\in {\cal T}_{2\kappa-1},  \deg(\varnothing)=1, \forall u\in T\setminus{\varnothing}, \deg(u)\in\{0,2\} \r\}$ be
the set of trees with $2\kappa-1$ edges, with binary branching points (except the root that has only one child). Denote by \[
\Delta_{n,M}=\{T\in {\cal T}_n,  \forall (u,v)\in \frak E(T), d(u,v)\in\sqrt{n}[M^{-1},M], T^b \in {\cal T}^B_{2\kappa}\}
.\]
A tree in $\Delta_{n,M}$ has its shape in ${\cal T}^B_{2\kappa}$ and all its spanned branches lengths in $\sqrt{n}[M^{-1},M]$.
\begin{lem}\label{qd} For any $`e>0$, there exists $M>0$ such that for $n$ large enough
\[`P_n(\Delta_{n,M} )\geq 1-`e.\]
\end{lem}
\proof  This is a consequence of $\bh_n\dd\bh=2\se/\sigma_{\mu}$ and the properties of $\se$~: $\se$ is a.s. non null on $(0,1)$, and  the local minima of $\se$ are a.s. all different (the continuum random tree is a.s. a binary tree). ~$\Box$\medskip

For any  $T\in \Delta_{n,M}$, $\Phi_T$ sends the nodes of $L(T)\setminus\{\varnothing\}$ on the leaves of $T^b$, and the nodes of $Z(T)$ on the internal nodes of $T^b$ (different from $\varnothing$), $\#Z(T)=\kappa-1$  (the branching nodes are distinct), $L(T)\cap Z(T)=\emptyset$, and $\Theta_T$ belongs to ${\cal D}^b={\cal D}_{2}\times {\cal D}_3^{\kappa -1}$ where ${\cal D}_{2}=\{(c,x), 1\leq x\leq c\}$ and  ${\cal D}_3=\{(c,x,y),  1\leq x< y \leq c\}$. Note that under $({\rm H}_1)$, ${\cal D}^b$ is a subset of $\cro{1,K}^{3\kappa-1}$.

\subsubsection{A second comparison result} 
\label{tg}
On the first hand, consider  $({\cal A}_n, {\cal H}_n,\Theta_n,\bt_n^b)$ the r.v.  $({\cal A}_\bt, {\cal H}_\bt,\Theta_\bt,\bt_n^b)$ when $\bt$ is $`P_n$-distributed. On the other hand, we define  $({\cal A}^{\star}_n,{\cal H}_n,\Theta^{\star}_n,\bt_n^b)$  as follows. Let  ${\cal H}_n=({\cal H}_n^{1},\dots,{\cal H}_n^{\#\frak{E}_n})$ and ${\cal A}_n^{\star}=({\cal A}_n^{i,\star})_{i=1,\dots,\#\frak{E}_n}$. Conditionally on ${\cal H}_n$,  the r.v. ${\cal A}_n^{i,\star}$s are independent with respective distribution $`Q_{{\cal H}_n^i}$.  The r.v.
$\Theta^{\star}_n=(\Theta^{\star}_n(i))_{i=1,\dots,\kappa}$ is independent of 
$({\cal H}_n,{\cal A}_n^{\star},\bt_n^b)$ and :
\be
`P\l(\Theta^{\star}_n(1)=(o,j^1)\r)&=&\mu_{o} \textrm{ for any } (o,j^1)\in {\cal D}_{2}\\
`P\l(\Theta^{\star}_n(i)=(o,j^1,j^{2})\r)&=&\tilde{\mu}_{o}:={\mu_{o}}/{\sigma^2} \textrm{ for any } (o_1,j^1,j^2)\in {\cal D}_3, i\geq 2.
\ee
Since the mean and the variance under $\mu$ are respectively 1 and $\sigma^2$, these formulas define indeed two distributions.
For any object $O$ and $l\in \mathbb{N}$, denote by $O^{[l]}=(O_1,\dots,O_l)$. Let \be\Gamma_{n,M}&=& \l\{\l(\sa^{[2\kappa-1]},\sx^{[2\kappa-1]},\theta^{[\kappa]},T^b\r),  \sx_i \in\sqrt{n}[M^{-1},M], \sa_i \in J_{\sx_i},T^b \in T_{2\kappa-1}^B \r\}\\
&&\cap \supp\l({\cal A}^{\star}_n,{\cal H}_n,\Theta^{\star}_n,\bt_n^b\r).\ee
The following Proposition generalizes to finite-dimensional distributions the Proposition \ref{ouch}.
\begin{pro} \label{yopyop}
$(i)$  For any $`e>0$, there exists $M$ such that 
\[`P\l(({\cal A}^\star_n,{\cal H}_n,\Theta_n^{\star},\bt^b_n)\in\Gamma_{n,M} \r)\geq 1-`e.\]
$(ii)$ \begin{equation}
\sup_{(\sa,\sx ,\theta,\tau_b)\in \Gamma_{n,M}} \l|\frac{`P\l(({\cal A}_n,{\cal H}_n,\Theta_n,\bt_n^b)=(\sa, \sx ,\theta,\tau^b)\r)}{`P\l(({\cal A}_n^{\star},{\cal H}_n,\Theta_n^\star,\bt_n^b)=(\sa, \sx ,\theta,\tau^b)\r)}
-1\r|\xrightarrow[n]{} 0.
\end{equation}
\end{pro}
\begin{com}\label{yopp} In general $`P_{({\cal A}_n,{\cal H}_{n},\Theta_n,\bt_n^b)}\nprec`P_{({\cal A}^\star_n,{\cal H}_n,\Theta^\star_n,\bt^b_n)}$ since $\supp(\Theta_n^{\star})$ is strictly included in $\supp(\Theta_n)$ when $\mu[3,+\infty)>0$ (the variable $\Theta^{\star}$ mimics the coding of binary branchings on $\bt_n^b$). In that case, moreover $`P_n(\bt_n^b\notin T_{2\kappa-1}^B)>0$ for $n$ large enough, and no control of $({\cal A}_n,{\cal H}_n,\Theta_n,\bt_n^b)$ is provided on $\complement T_{2\kappa-1}^B$. 
Notice that the condition $`P_{X_n}\prec `P_{Y_n}$ in Definition \ref{def1} may be replaced by the following weaker condition sufficient to keep the conclusion of Lemma \ref{dr}~: 
for any $`e>0$, for $n$ large enough, there exists a measurable set $A_n^{`e}$ such that $`P_{Y_n}(A_n^{`e})\geq 1-`e$ and a function $f_n^{`e}: A_n^{`e}\mapsto `R$ satisfying $`P_{X_n}=f_n^{`e} `P_{Y_n}$ on $A_n^{`e}$ and $\sup|f_n^{`e}-1|<`e$. Here, this is the case on  $\Gamma_{n,M}$. 
\end{com}

\proof $(i)$ is a consequence of Lemmas \ref{qd} and \ref{degdeg}. \\
For $(ii)$, let  $(\sa,\sx,\theta,\tau^b)\in \Gamma_{n,M}$ for $\theta=\l((o_1,j^1_1),(o_2,j_2^1,j_2^2),\dots,(o_{\kappa-1},j^1_{\kappa-1},j^2_{\kappa-1})\r)$. By classical properties of GW trees $`P\l(({\cal A}_n,{\cal H}_n,\Theta_n,\bt_n^b)=(\sa, \sx ,\theta,\tau^b)\r)$
\ben\label{cost}
&=& \frac{\l(\dis\prod_{i=1}^{2\kappa-1} `Q_{{\sx_i}}(\sa_i)\r) \l(\dis\prod_{i=1}^\kappa \mu_{o_i}\r)`P\l(|{\bf f}^{(l)}_{\#\Sub(\tau^b,\sa,\theta)}|=F_l(\sx ,\tau^b),0\leq l\leq \kappa+1\r) }{`P(|\bt|=n)},
\een
where the ${\bf f}^{(l)}$'s are independent forests.
On the other hand, $`P\l(({\cal A}^{\star}_n,{\cal H}_n,\Theta_n^\star,\bt_n^b)=(\sa, \sx ,\theta, \tau^b)\r)$
\be
&=& `P\l(({\cal A}^{\star}_n,\Theta_n^\star)=(\sa,\theta)\, \b|\, \b({\cal H}_n,\bt_n^b\b)= (\sx , \tau^b\b)\r)`P\l(\b({\cal H}_n,\bt_n^b\b)= (\sx , \tau^b\b)\r)\\
&=& \l(\dis\prod_{i=1}^{2\kappa-1} `Q_{{\sx_i}}(\sa_i)\r) \l(\prod_{i=1}^\kappa \tilde{\mu}_{o_i}\r)`P\l(\b({\cal H}_n,\bt_n^b\b)= (\sx , \tau^b\b)\r);
\ee
summing formula \eref{cost} on all possible values of the $\sa_i$'s and the $\theta$'s leads to
\begin{equation}\label{mul}
`P\l(({\cal H}_n,\bt_n^b)= (\sx , \tau^b)\r)=\frac{\sigma^{2(\kappa-1)}\,`P\l(|{\bf f}^{(l)}_{\#\Sub(\tau^b,{\sf m},\tilde{\theta})}|=F_l(\sx ,\tau^b),0\leq l\leq\kappa+1\r)}{`P(|\bt|=n)},
\end{equation}
where ${\sf m}=({\sf m}_i)_{i=1,\dots,\kappa}$ is a vector of $\kappa$ multinomial independent r.v. (the parameters of ${\sf m}_i$ are $\sx _i$ and ${\bf p}$), and where $\tilde\theta=(\tilde\theta(i))_{i\in\cro{1,\kappa}}\sur{=}{(d)}\Theta_n^{\star}$ and is independent of ${\sf m}$.
Hence,  for $(\sa,\sx,\theta,\tau^b)\in \Gamma_{n,M}$,
\ben
\label{rap}
\frac{`P\l(({\cal A}_n,{\cal H}_n,\Theta_n,\bt_n^b)=(\sa, \sx ,\theta,\tau^b)\r)}{`P\l(({\cal A}_n^{\star},{\cal H}_n,\Theta_n^\star,\bt_n^b)=(\sa, \sx ,\theta,\tau^b)\r)}
&=&\frac{`P\l(|{\bf f}^{(l)}_{\#\Sub(\tau^b,\sa,\theta)}|=F_l(\sx ,\tau^b),0\leq l\leq\kappa+1 \r)}{`P\l(|{\bf f}^{(l)}_{\#\Sub(\tau^b,{\sf m},\tilde{\theta})}|=F_l(\sx ,\tau^b),0\leq l\leq\kappa+1\r)}.
\een

It is easy to check that for any $(\sa,\sx,\tau^b ,\theta)$ in $\Gamma_{n,M}$,  any $l$, for $n$ large enough
 \[|F_l(\sx ,\tau^b)-n(s_{l+1}-{s_l})|\leq n^{2/3},~~~|\#\Sub_l(\tau^b,\sa,\theta)-\frac{\sigma^2}{2}d(u_l,u_{l+1})|\leq n^{5/12},\]
since $(1/2)^{2/3}<5/12$. This allows to approximate on one hand $F_l(\sx ,\tau^b)$ by $n(s_{l+1}-{s_l})$, and on the other hand $\#\Sub_l(\tau^b,\sa,\theta)$ by $\frac{\sigma^2}{2}d(u_l,u_{l+1})$ on $\Gamma_{n,M}$ (since $n^{5/12}=o(n^{1/2})$, the order of $d(u_l,u_{l+1})$).
So, using Otter and the central local limit theorem and also a decomposition of the denominator along $\{{\sf m} \in \prod J_{\sx_i}\}$ or in its complements (as in the proof of Proposition \ref{RN}), we get 
\begin{equation}\label{erzp}
\sup_{(\sa,\sx ,\theta,\tau_b)\in \Gamma_{n,M}} \l|\frac{`P_n\l(|{\bf f}^{(l)}_{\#\Sub(\tau^b,\sa,\theta)}|=F_l(\sx ,\tau^b),0\leq l\leq \kappa+1\r)}{`P_n\l(|{\bf f}^{(l)}_{\#\Sub(\tau^b,{\sf m},\tilde{\theta})}|=F_l(\sx ,\tau^b),0\leq l\leq\kappa+1\r)}-1\r|\xrightarrow[]{n} 0.~~\Box
\end{equation}

\subsubsection{Proof of the convergence of the finite-dimensional distribution in Theorem \ref{zozo}}
We now show that  Proposition \ref{yopyop} implies the convergence of the finite-dimensional distributions in Theorem \ref{zozo}. The proof is similar to the one of Corollary \ref{opieds}. \par
Thanks to the Skohorod representation theorem \cite[Theorem 3.30]{KAL}, there exists a probability space  $\Omega$ on which the convergence of $\bh_n$ to $\bh$ is a.s..
On $\Omega$,  the vector 
\[V_n=(\bh_n(s_1),\check{\bh}_n(s_1,s_2),\bh_n(s_2),\check{\bh}_n(s_2,s_3),\dots,\bh_n(s_\kappa)),\] which determines $\bt_n^b$ as well as the length of the spanned branches,  converges a.s. to 
\[V_{\infty}=(\bh(s_1),\check{\bh}(s_1,s_2),\bh(s_2),\check{\bh}(s_2,s_3),\dots,\bh(s_\kappa)),\] which determines $\tau_{\bf s}$ the subtree of the continuum random tree $\tau_{\infty}$, with contour process $\bh$, spanned by the root and the nodes visited at times $s_1,\dots,s_{\kappa}$ (see Aldous  \cite{ALD,ALD3}).
With probability 1, the coordinates of $V_{\infty}$ are distinct and non zero, and then $\tau_{\bf s}$ has its shape $\tau_{\bf s}^b$ in ${\cal T}_{2\kappa-1}^B$. Let ${\cal H}_\infty=({\cal H}_{\infty,i})_{i\in\cro{1,2\kappa-1}}$ be the lengths of the (sorted) spanned branches in $\tau_{\infty}$.
On $\Omega$, $\bt_n^b\as \tau_{\bf s}^b$, and for $M$ large enough, ($M$ depending on $\tau_\infty)$, and $n$ large enough,  $\bt_n\in \Gamma_{n,M}$. 

Denote by $(A_n^i)_{i\in\cro{1,2\kappa-1}}$ the (sorted) corresponding content of the spanned branches of $\bt_n$, and by $({\cal H}_n^i)_{i\in\cro{1,2\kappa-1}}=(|A_n^i|)_{i\in\cro{1,2\kappa-1}}$ their lengths.  The normalized contents are then given by 
\[\bg^{(n)}_{i,k,j}=n^{-1/4}\l((A_{n}^i)_{k,j}-\mu_k|{\cal H}_n^i|\r).\]
A consequence of Proposition \ref{yopyop}, is that 
\[\l((\bg^{(n)}_i)_{i\in\cro{1,2\kappa-1}},{\cal H}_n/\sqrt{n}\r)\di\l(({\cal G}^{(i)}(n,{\cal H}_n^i))_{i\in\cro{1,2\kappa-1}},{\cal H}_n/\sqrt{n}\r)\to 1,\]
in the sense of Comment \ref{yopp} (which slightly modifies Definition \ref{def1}) where the r.v. ${\cal G}^{(i)}(n,{\cal H}_n^i)$'s are independent, and conditionally on ${\cal H}_n^i=l$,  $({\cal G}^{(i)}(n,{\cal H}_n^i))\sur{=}{(d)}{\cal G}(n,l)$.  On $\Omega$,  ${\cal H}_n/\sqrt{n}\as{\cal H}_\infty$, and then by Proposition \ref{oo},  $({\cal G}^{(i)}(n,{\cal H}_n^i))_{i\in\cro{1,2\kappa-1}}$ converges in distribution to a centered Gaussian vector $({\cal G}^{(i)}_{\infty})_{i\in\cro{1,2\kappa-1}}$ with independent coordinates, where ${\cal G}^{(i)}_{\infty}$ has variance ${\cal H}_{\infty,i}$. This implies the convergence of the finite-dimensional distributions  in Theorem \ref{zozo}. $~\Box$

\subsection{Tightness in Theorem \ref{zozo}}
We only prove the tightness of the family $(\bG^n)$, since one already knows that $(\bh_n)$ is tight (since $\bh_n\dd \bh$). In this section, we assume $({\rm H}_1)$ and $({\rm H}_2)$.\par
We collect in the set  $\Omega_n^{\alpha,\delta,\gamma,\rho}$,  the trees with $n$ edges  having some suitable properties~: 
\be
\Omega_n^{\alpha,\delta,\gamma,\rho}&=&\Big\{T \in {\cal T}_n, \forall t,s \in[0,1], |\bh_n(s)-\bh_n(t)|\leq \delta|t-s|^{\alpha}, \max_{l} \b||u(l+1)|-|u(l)|\b|\leq \rho \log n,\\
&&|u(n)|<\rho\log n, \forall (k,j) \in I_K,  l\in(0,|u|], \left|A_{u,l,k,j}-\mu_k l\right|\leq \gamma\sqrt{l\log n} \Big\}
\ee
\begin{lem}\label{gros}
For any $`e>0$,  $\alpha<1/2$,  there exists $\delta>0$, $\gamma>0$, $\rho>0$, s.t. $`P(\Omega_n^{\alpha,\delta,\gamma,\rho})\geq 1-`e$. 
\end{lem}
According to Lemmas \ref{degdeg} and \ref{aukj}, and Remark \ref{rzer}, only the condition on the H\"olderienity of $H$ has to be checked. This is postponed at the end of the paper.\medskip

Let $`e>0$ be fixed. Set  $\alpha=2/5$ and choose $\delta>0$,  $\gamma>0$, $\rho>0$ s.t.  $`P(\Omega_n^{\alpha,\delta,\gamma,\rho})\geq 1-`e$ for $n$ large enough. For these choices, write $\Omega_{`e}$ instead of $\Omega_n^{\alpha,\delta,\gamma,\rho}$.

We will establish the following Proposition. 
\begin{pro}\label{oye}
There exists $a>0$,  $\beta>0$, $c>0$  s.t. for $n$ large enough, 
\begin{equation}\label{oye2}
`E\l(\|\bG_s^n-\bG_t^n\|_1^\beta`I_{\Omega_{`e}}\r)\leq c\,|t-s|^{1+a} ~~~\textrm{ for any }s,t\in[0,1].
\end{equation}
\end{pro}
This implies that the $(1+a)/\beta$-H\"older norm of the family $(\bG^n)$ is tight, and then that  $(\bG^n)$ is tight in $C([0,1])^{\#I_K}$ (recall that $\bG_0^n$ is the null vector of $`R^{\#I_K}$). \par
We first point out that using \eref{mau}, we get that for any $a>0$ there exists $\beta>0$ such that for $n$ large enough
$ `E( \|\bGn(s_n)-\bGn(s)\|_1^\beta`I_{\Omega_n^{\rho}}) \leq  c(s-s_n)^{1+a}.$ Hence, we can  restrict ourself to prove \eref{oye2} only for $s$ and $t$ such that $ns$ and $nt$ are integer (this is classical).
From now on, we assume that $s,t$ are in $[0,1]_n:=[0,1]\cap \mathbb{N}/n$, and $s\neq t$.\par
We set  $u_1=u(\floor{ns}),u_2=u(\floor{nt})$, $\check{u}_{1,2}$ their deepest common ancestor and  $D_n(s,t)=d(u_1,u_2)$. There exists $\delta'>0$, such that for $T \in\Omega_{`e}$, any $s,t\in[0,1]_n$, $s\neq t$, 
\[D_n(s,t)\leq 2+H_n(nt)+H_n(ns)-2\min_{k\in [ns,nt]}H_n(k)\leq 2\sqrt{n}\,\delta|t-s|^{2/5}+2\leq \delta'\sqrt{n}|t-s|^{2/5}. \]

\begin{lem}\label{yuplm}
For  any  $\alpha'>0,a>0$, there exists $\beta>0,c>0$ s.t. for any $s,t\in[0,1]$ such that $|s-t|\leq (\log n)^{-3}$, for $n$ large enough,
\ben
`E\l(\|\bG_s^n-\bG_t^n\|_1^\beta\,\,`I_{\Omega_{`e}}\r)&\leq& c |t-s|^{1+a}.
\een
\end{lem}
\proof  Let $s,t\in[0,1]_n$, $s\neq t$. 
We use a deterministic bound valid for all trees $T$ in $\Omega_{`e}$. 
Let $(k,j)\in I_K$ fixed. 
As in the proof of Proposition \ref{qe}, it suffices to show that 
\begin{equation}\label{azeaze}
{n^{-\beta/4}\l|A_{u_1,k,j}-\mu_k|u_1|-A_{u_2,k,j}+\mu_k|u_2|\r|^\beta}\leq c|s-t|^{1+a}
\end{equation}
Passing via $\check{u}_{1,2}$, the left hand side of \eref{azeaze} is smaller than
\[c_1n^{-\beta/4}\l(|A_{u_1,h_1,k,j}-\mu_k|h_1||^\beta+|A_{u_2,h_2k,j}-\mu_k|h_2||^\beta +2^\beta\r)\]
where $h_1:=d(u_1,\check{u}_{1,2})-1$ and $h_2:=d(u_2,\check{u}_{1,2})-1$ (the contribution of  $\check{u}_{1,2}$ is bounded by the term 2). Using that $\left|A_{u,l,k,j}-\mu_k l\right|\leq \gamma\sqrt{l\log n}$ for any $l$ and $l\leq D_n(s,t)\leq \delta'{n^{1/2}}|t-s|^{2/5}$, we find 
\[{n^{-\beta/4}\l|A_{u_1,k,j}-\mu_k|u_1|\r|^{\beta}+\l|A_{u_2,k,j}-\mu_k|u_2|\r|^\beta}\leq c_2(t-s)^{\beta/5}(\log n)^{\beta/2}\]
and since $|t-s|\leq (\log n)^{-3}$, $|t-s|^{\beta/6}(\log n)^{\beta/2}\leq 1$ and then 
$c_2(t-s)^{\beta/5}(\log n)^{\beta/2}$ is smaller than $|t-s|^{1+a}$ for $\beta$ and $n$ large enough.~$\Box$\medskip

\begin{lem} For any $a>0$, there exists $\beta>0$, $c>0$ s.t.
for any $t\in[0,1]$, for any $n$ large enough,
\begin{equation}\label{mauve}
`E\l(\|G_t^n\|_1^\beta\,\1_{\Omega_{`e}}\r)\leq c \, t^{1+a}.\end{equation}
\end{lem}
\proof Consider first the case $t=1$. In $\Omega_\varepsilon$, we have $|u(n)|\leq \rho \log n$ and then
\[`E(\|\bG_1^n\|_1^{\beta}\1_{\Omega_{`e}})
\leq c_1(\rho \log n)^{\beta} n^{-\beta/2}\]
and this is smaller than $c 1^{1+a}$ for any $a>0$, $c>0$, $\beta>0$ for $n$ large enough. By the previous Lemma and a simple computation (using that $\max ||u(l+1)|-|u(l)||\leq \rho\log n$) one sees that \eref{mauve} is true if $t\notin V_n$ where
\[V_n:=[(\log n)^{-3},1-(\log n)^{-3}].\]
Assume now that $t\in V_n$.
In $\Omega_{`e}$, the H\"older property of $\bh_n$ and the inequality $|u(n)|\leq \rho \log n$, implies that for $t\in V_n$,
\begin{equation}\label{vivo}
|u(\floor{nt})|\leq \bar{L_n}(t):=c_2 n^{1/2}[t\wedge(1-t)]^{\alpha}.
\end{equation}

For any real number $a$, we denote by $a.\mu$ the vector $(a\mu_k)_{(k,j)\in I_K}$.
Using \eref{ouch} and \eref{ptn}, there exists $c_3>0$ such that for $t\in V_n$ and $n$ large enough, $`E\l(\|G_t^n\|_1^\beta\,\1_{\Omega_{`e}}\r)\leq $
\[
c_3 \sum_{h\leq\bar{L_n}(t)}\sum_{{\sf a}\in \mathbb{N}^I_{[h]}}`Q_{h}(\sa)\,\frac{\|\sa -h.\mu\|_1^\beta}{n^{\beta/4-3/2}}`P\l(|{\bf f}_{N_1(\sa)}|=\floor{nt}-h,\,|{\bf f}\,'_{1+N_2(\sa)}|=n+1-\floor{nt}\r)\]
and by Otter than
\[c_4 \sum_{h\leq\bar{L_n}(t)}\sum_{{\sf a}\in \mathbb{N}^I_{[h]}}`Q_{h}(\sa)\frac{\|\sa -h.\mu\|_1^\beta}{n^{\beta/4-3/2}} \frac{N_1(\sa)(1+N_2(\sa))`P(W_{\floor{nt}-h}=N_1(\sa))`P(W_{n-\floor{nt}+1}=1+N_2(\sa))}{(\floor{nt}-h)(n-\floor{nt}+1)},\]
where $(W_k)$ is the random walk described in the beginning of Section \ref{zrf}.
In order to bound these two last probabilities, we use a classical concentration property valid for any non-degenerate random walk $(W_k)_k$ (trivial consequence of Petrov \cite[Theo. 2.22 p.76]{PET})~: there exists a constant $c_5$ such that for any  $n\geq 0$,
\begin{equation}\label{zare}
\sup_y `P(W_n= y)\leq c_5/{\sqrt{n}}.
\end{equation}
Now, for any $\sa\in `N^I_{[h]}$, $N_1(\sa)$ and $N_2(\sa)$ are smaller than $Kh$,
and  for any $h\leq \bar{L_n}(t)$, $t\in V_n$ and $n$ large enough, $nt-h\geq nt/2$.
 We then get 
\be
`E(\|G_t^n\|_1^\beta\1_{\Omega_{`e}})&\leq& 
c_6\sum_{h\leq \bar{L_n}(t)}\sum_{{\sf a}\in \mathbb{N}^I_{[h]}}\frac{ `Q_{h}(\sa)\|\sa -h.\mu\|_1^\beta\,\,h^2}{n^{\beta/4-3/2}(\floor{nt}-h)^{3/2}(n-\floor{nt}+1)^{3/2}}.
\ee
Using Proposition \ref{qe}, we obtain that for any $t\in V_n$,
\[`E(\|G_t^n\|_1^\beta\1_{\Omega_{`e}})\leq 
c_7 \frac{\l(\bar{L_n}(t)\r)^{\beta/2+3}}{n^{\beta/4+3/2}(t(1-t))^{3/2}}\leq c_8 \frac{(t\wedge(1-t))^{\beta/2+3}}{(t(1-t))^{3/2}}.~\Box \]
\begin{remi}The last formula implies that for any $a>0$, there exists $\beta>0$, $c>0$ s.t.
for any $t\in V_n$, for any $n$ large enough,
\begin{equation}\label{concler}
`E(\|G_t^n\|_1^\beta\1_{\Omega_{`e}})\leq c (t\wedge(1-t))^{1+a}.
\end{equation}
This allows to prove a part of Proposition \ref{oye}~:
since $
`E(\|G_t^n-G_s^n\|_1^\beta\1_{\Omega_{`e}})\leq c`E(\1_{\Omega_{`e}}(\|G_t^n\|_1^\beta+\|G_s^n\|_1^\beta))$
when $s,t\in V_n$ and $s\leq t$,\\
-- if $s\leq t-s$ (in this case $t\leq 2(t-s)$) then $`E(\|G_t^n-G_s^n\|_1^\beta\1_{\Omega_{`e}})\leq c(t-s)^{1+a}$, \\
-- if $1-t\leq t-s$ (in this case $1-s\leq 2(t-s)$) then 
\[`E(\|G_t^n-G_s^n\|_1^\beta\1_{\Omega_{`e}})\leq c((1-t)^{1+a}+(1-s)^{1+a})\leq c_2(t-s)^{1+a}.\]
\end{remi}
Thanks to this remarks, only the case  $s,t\in V_n$, $s\leq t$, and 
\begin{equation}\label{gfd}
[s\wedge(1-s)]\geq t-s \textrm{ and }[t\wedge(1-t)]\geq t-s
\end{equation}
remains to be checked. So assume that $s$ and $t$ satisfies these constraints. 

Consider 
${\cal A}_n=({\cal A}_n^1,{\cal A}_n^2,{\cal A}_n^3)=(A_{( \check u_{1,2}, u_{1})},A_{( \check u_{1,2}, u_{2})},A_{(\varnothing, \check u_{1,2})})$ the contents of the ``three'' spanned branches in $\bt_{{\bf s}^2}$ (some of these spanned branches may be empty). 
We have
\ben\label{prem}
`E\l(\|\bG_s^n-\bG_t^n\|_1^\beta`I_{\Omega_{`e}}\r)\leq\sum_{h_1,h_2,h_3}\sum_{\sa_1,\sa_2,\sa_3}\frac{`P_n({\cal A}_n^i=\sa_i,i=1,2,3)\l[\|\sa_1-h_1.\mu\|_1^\beta+\|\sa_2-h_2.\mu\|_1^\beta\r]}{n^{\beta/4}}
\een
where the first sum is taken on $h_1+h_3\leq \bar{L_n}(s)$, $h_2+h_3\leq \bar{L_n}(t)$, $h_1+h_2\leq \bar{D_n}(s,t):=\delta'n^{1/2}|t-s|^{\alpha}$ where $\bar{L_n}(x)$ is given in \eref{vivo}. 
By Section \ref{subs}, and the Otter formula, $h_1,h_2,h_3,\sa_1,\sa_2,\sa_3$ fixed,
\begin{equation}\label{eto}
`P_n({\cal A}_n^i=\sa_i,i=1,2,3)\leq cn^{3/2}\sup_{\theta}\prod_{i=1}^3`Q_{h_i}(\sa_i)\frac{S_i(\theta)`P(W_{F_i}=S_i(\theta))}{F_i}
\end{equation}
where the supremum is taken on $\theta=(\theta_1,\theta_2,\theta_3)\in \cro{0,K}^3$, and where
$F_1=ns+1-|u(\floor{ns})|-1$, $F_2=n(t-s)+1-(|u(\floor{nt})|-|\check{u}_{1,2}|)$, $F_3=n(1-t)+1$, $S_1=N_1(\sa_3)+N_1(\sa_1)+\theta_1$, $S_2=N_2(\sa_1)+N_1(\sa_2)+\theta_2$, $S_3=N_2(\sa_2)+N_2(\sa_3)+\theta_3$.  

We plug this bounds in \eref{prem}, and bound the left hand side using the following ingredient:\\
-- the probabilities in \eref{eto} involving the random walks are bounded using \eref{zare}. \\
-- for $\sa\in `N^I_{[h]}$, $N_i(\sa)\leq Kh$ and then for a constant $c>0$,
\be 
S_1 &\leq & K|u(\floor{ns})|+\theta_1\leq c\bar{L_n}(s),\\
S_2 &\leq &c |\bar{D_n}(s,t)|,\\
S_3&\leq & K|u(\floor{nt})|+\theta_3\leq c\bar{L_n}(t).
\ee
The denominator are bounded using $|t-s|\geq (\log n)^{-3},[t\wedge(1-t)]\geq (\log n)^{-3},[s\wedge(1-s)]\geq (\log n)^{-3}$, and then for $n$ large enough,
\[F_1 \geq  ns/2,\quad F_2 \geq  n(t-1)/2,\quad F_3 \geq  n(1-t)/2.\]
Finally we get that the left hand side of \eref{prem} is smaller than
\be
c\frac{\bar{L_n}(s)\bar{L_n}(t)\bar{D_n}(s,t)\sum_{h_1,h_2,h_3}\sum_{\sa_1,\sa_2,\sa_3}\prod_{i=1}^3`Q_{h_i}(\sa_i)\l[\|\sa_1-h_1.\mu\|_1^\beta+\|\sa_2-h_2.\mu\|_1^\beta\r] }{n^{\beta/4-3/2}\l[n^{3}(s\wedge(1-s))(t\wedge(1-t))(t-s)\r]^{3/2}}
\ee
The double sum is smaller than
\[\sum_{h_1,h_2,h_3} h_1^{\beta/2}+h_2^{\beta/2}\leq (\bar{D_n}(s,t))^{\beta/2+2} \bar{L_n}(s)\]
this last factor $\bar{L_n}(s)$ being a bound of $h_3$. 
Finally,
\[`E\l(\|\bG_s^n-\bG_t^n\|_1^\beta`I_{\Omega_{`e}}\r)\leq cn^{3/2-\beta/4}\frac{(\bar{L_n}(s))^2\bar{L_n}(t)(\bar{D_n}(s,t))^{\beta/2+3}}{\l[n^3 (s\wedge(1-s))(t\wedge(1-t))(t-s)\r]^{3/2}}\]
By \eref{gfd}, it suffices to take $\beta$ large enough.~$\Box$

\subsection{Proof of Theorem \ref{serpdet}}
Consider the representation of $\ell(u)$ given in \eref{ellu}. For any $s$ such that $ns$ is an integer,
\begin{equation}\label{ellup}
\br_n(s)=\br_n^{(1)}(s)+\br_n^{(2)}(s)
\end{equation}
where
\be
\br_n^{(1)}(s)&=&n^{-1/4}\sum_{(k,j)\in I_K}\sum_{l=1}^{A_{u(ns),k,j}} \b(Y_{k,j}^{(l)}-m_{k,j}\b),\\
\br_n^{(2)}(s)&=&n^{-1/4}\sum_{(k,j)\in I_K}\l(A_{u(ns),k,j}-\mu_k |u(ns)|\r)m_{k,j}=<\bG^{(n)}(s), \ov{m}>,
\ee
where $\ov{m}=(m_{k,j})_{(k,j)\in I_K}$ and $<a,b>=\sum_{(k,j)\in I_K}a_{k,j}b_{k,j}$. For $s$ in $[i/n,(i+1)/,n]$, $\br_n^{(1)}(s)$ and $\br_n^{(2)}(s)$ are defined by linear interpolation.
Since $\bh_n\dd\bh$ in $C[0,1]$, by the Skohorod representation theorem \cite[Theorem 3.30]{KAL}, there exists a probability space $\Omega$ on which this convergence is a.s.. On this space  by Theorem \ref{zozo},  $\bG^{(n)}$  converges in distribution in $C([0,1])^{\#I_K}$ to $\bG^{\bh}$, where $\bG^{\bh}$ has the distribution of $\bG$ knowing $\bh$.
Now, since   the application 
\[\app{\Psi_{\ov m}}{C([0,1])^{\#I_K}}{(C[0,1])}{\b(s\mapsto g(s)\b)}{\b(s\mapsto<g(s), \ov{m}>\b)}\]
is continuous,  on $\Omega$ we have
\begin{equation}
<\bG^{(n)}, \ov{m}>~\dd~ \br^{(2)}:=<\bG^{\bh}, \ov{m}>
\end{equation}
in $C([0,1])$. On $\Omega$, $\br^{(2)}$ is a centered Gaussian process with covariance function 
\[\cov\b(\br^{(2)}(s),\br^{(2)}(t)\b)=\check\bh(s,t)\sum_{(k,j)\in I_K}\sum_{(k',j')\in I_K} (-\mu_k\mu_{k'}+\mu_k\1_{(k,j)=(k',j')})m_{k,j}m'_{k,j}.\]
On the other hand, $\br^{(1)}_n$ is the standard head of a discrete snake associated with independent displacements. As shown in \cite{GM}, under $({\rm H}_1)$ and $({\rm H}_2)$,
\begin{equation}
\br_n^{(1)}\dd\br^{(1)}
\end{equation}
in $C([0,1],`R)$ where $\br^{(1)}$ is a centered Gaussian process with covariance function   
\[\cov\b(\br^{(1)}(s),\br^{(1)}(t)\b)=\check{\bh}(s,t)\sum_{(k,j)\in I_K}\mu_k \sigma^2_{k,j}.\]

We shall prove that, given $\bh$, the finite-dimensional distributions of $\br^{(1)}$ and $\br^{(2)}$ are independent. We establish the ``asymptotic independence'' between the two processes $\br_n^{(1)}$ and $\br_n^{(2)}$ knowing $\bh$.  The arguments are quite straightforward; we just explicit the uni-dimensional case.
Let 
\[{\cal T}_n^{\nu}=\l\{T\in{\cal T}_n,  \forall (k,j)\in I_K,  u\in{T}, \big|A_{u,k,j}-\mu_k |u|\big|\leq  n^{1/4+\nu}\right\}.\]
According to Lemma 9 in \cite{GM} (it is also a consequence of Lemma \ref{gros}), for any $\nu>0$, $`e>0$, if $n$ is large enough $`P_n( {\cal T}_n^{\nu})\geq 1-`e$. Let $s\in[0,1]$ (such that $ns$ is an integer), one may compare 
\be
\br_n'(s)={n^{-1/4}\sum_{(k,j)\in I_K}\sum_{l=1}^{\floor{\mu_k |u(\floor{ns})|- n^{1/4+\nu}} } \b(Y_{k,j}^{(l)}-m_{k,j}}\b)
\ee
where the same r.v.  $Y_{k,j}^{(l)}$ are involved in both $\br_n'$ and $\br^{(1)}_n$, with $\br_n^{(1)}(s)$. Since knowing $|u(\floor{ns})|$, $\br_n'(s)$ is clearly independent of $\br_n^{(2)}(s)$, the proof of $|\br_n^{(1)}(s)-\br'_n(s)|\proba 0$ will prove our claim (in the uni-dimensional case). 
We have
\be
`P_n(|\br_n'(s)-\br_n^{(1)}|\geq x)&\leq &`P(|\br_n'(s)-\br_n^{(1)}|\geq x, 
{\cal T}_n^\nu)+`P_n({\cal T}_n\setminus{\cal T}_n^\nu).
\ee
The  last term goes to 0 for any  $\nu>0$. The Rosenthal inequality \cite[Theorem 2.11]{PET} asserts that if $(X_k)_k$ is a sequence of centered r.v.  and $q\geq 2$, then 
\begin{equation}\label{rose}
`E\Big(|\sum_{i=1}^n X_i|^{q}\Big)\leq c(q)\Big(\sum_{i=1}^n `E(|X_i|^q) + \Big(\sum_{i=1}^n \var(X_i)\Big)^{{q}/2}\Big)\end{equation}
where $c(q)$ is a positive constant depending only on $q$. For $p$ satisfying $({\rm H}_2)$, we have  
\[
`P(|\br_n'(s)-\br_n^{(1)}|\geq x, 
{\cal T}_n^\nu)\leq `E\l(x^{-p}|\br_n'(s)-\br_n^{(1)}|^{p}
\1_{{\cal T}_n^\nu}\r)
\]
Conditioning at first by the $A_{(u(ns))}$, and using \eref{rose}, we get $`P(|\br_n'(s)-\br_n^{(1)}|\geq x, {\cal T}_n^\nu)\leq$
\[\frac{x^{-p}c(p)}{n^{{p}/4}}\Big(\sum_{(k,j)\in I_K} 2 n^{1/4+\nu} `E(|Y_{k,j}-m_{k,j}|^{p}) +\Big(\sum_{(k,j)\in I_K} 2 n^{1/4+\nu} \sigma_{k,j}^2)\Big)^{p/2}\Big)\]
and then for $\nu<1/4$, for any $x>0$ the bound goes to 0.\par
Hence $\br$ is a centered Gaussian process with covariance function sum of the ones of $\br^{(1)}$ and $\br^{(2)}$. Using that $\sum_{(k,j)}\mu_km_{k,j}={\bf m}=0$, we get $\cov(\br(s),\br(t))=\check{\bh}(s,t)\sum_{(k,j)}\mu_k`E(Y_{k,j}^2).~~~\Box.$

\subsection{Proof of Lemma \ref{gros}} We say that a sequence of processes $(u_n)$ defined on $[0,1]$ is uniformly H\"older continuous in probability with exponent $\alpha$ ($\alpha$-UHCP) if 
for every $`e>0$ there exists a real number $C_{`e}$ such that for every $n$, 
\[`P(|u_n(s)-u_n(t)|\leq C_{`e}|t-s|^{\alpha}
\textrm{ for all }s,t\in[0,1])\geq 1-`e.\] 

Only the following condition on $(\bh_n)$ remains to be checked~:
 \begin{lem} $({\bh_n})$  is $\alpha$-UHCP for any $\alpha<1/2$.
\end{lem}
\bf Proof~: \rm For any tree $T\in{\cal T}$, $H_T$ is a simple function of $\widehat{H}_T$~: let $m_T(0)=0$, and for any $i\geq 1$, $m_T(i)=\min\{j, j>m_T(i-1), \widehat{H}_T(j)>\widehat{H}_T(j-1)\}$, then
$H_T(k)=\widehat{H}_T(m_T(k)).$ In fact, $m_T(k)=\inf\{j, F_T(j)=u(k)\}$. One may check inductively on $k$ that,
\begin{equation}\label{mj}
m_T(k)+H_T(k)=2k \textrm{ for any }k\geq 0.
\end{equation}

Assume that $|\widehat{\bh_n}(s)-\widehat{\bh_n}(t)|\leq c_1|t-s|^{\alpha}$ and $|{\bh_n}(s)-{\bh_n}(t)|\leq c_2|t-s|^{\beta}$, for some $\alpha,c_1,c_2\geq 0, \beta\in[0,1/2]$ for any $s,t\in[0,1]$. 
Then for $s,t \in [0,1]_n$,
\be
|\bh_n(s)-\bh_n(t)|&=&\l|\widehat{\bh}_n(m(ns)/(2n))-\widehat{\bh}_n(m(nt)/(2n))\r|
                   \leq c_1\l|\frac{m(ns)-m(nt)}{2n}\r|^{\alpha}\\
&\leq&  c_1\l|{s-t}   + \frac{H(nt)-H(ns)}{2n}\r|^{\alpha}
\leq  c_1\l|{s-t}   + c_2 \frac{|t-s|^{\beta}}{\sqrt{n}}\r|^{\alpha}\\
&\leq&  c_1\l|{s-t}   + c_2 |t-s|^{\beta+1/2}\r|^{\alpha}.
\ee
Since when $s,t\in [0,1]_n$ and $\beta<1/2$, we have $|s-t|\leq  |t-s|^{\beta+1/2}$. Hence, for any $s,t\in[0,1]$, by interpolation, we get that $|\bh_n(s)-\bh_n(t)|\leq c_3 |t-s|^{\alpha(\beta+1/2)}$ for a certain constant $c_3$. Hence, if $\widehat{\bh}_n$ is $\alpha$-UHCP and $\bh_n$ is $\beta$-UHCP then $\bh_n$ is $f_{\alpha}(\beta):={\alpha(\beta+1/2)}$-UHCP.\par

By Gittenberger \cite{G1}, for all $s$,
$t$, $`e>0$  
\begin{equation*}
  `P(|\widehat{\bh_n}(s)-\widehat{\bh_n}(t)|\ge`e) \le C_1|s-t|^{-1}\exp\l(-C_2`e|s-t|^{-1/2}\r),
\end{equation*}
which  gives, for any $p>0$,  $`E\l(\l|\widehat{\bh_n}(s)-\widehat{\bh_n}(t)\r|^p\r) \le C(p)|s-t|^{p/2-1}.$ Taking $p$ large enough, this ensures that for every $\alpha<1/2$, the family $(\widehat{\bh_n})$ is $\alpha-$UHCP. \par
Fix one such $\alpha$. Since $\bh_n$ is 0-UHCP it is also $f_{\alpha}(0)$-UHCP, and by successive iterations $f_{\alpha}^m(0)$-UHCP. Using that $f_{\alpha}$ is increasing and has $\beta_{\alpha}=\alpha/(2(1-\alpha))$ as finite fix point, a finite number of iterations shows that $\bh_n$ is $\beta$-UHCP for every $\beta<\beta_{\alpha}$. Finally, since $\lim_{\alpha\to1/2}\beta_\alpha=1/2$, for any $\beta<1/2$, we can choose $\alpha<1/2$ such that $\beta<\beta_{\alpha}$.~$\Box$

\subsection{Note on the comparison between height processes and contour processes}
\label{note}
Let $T$ be a tree in ${\cal T}_n$, $F_T$ its depth first traversal,  and $(u(k))_{k\in\cro{0,n}}$, the sorted list of its nodes. 
Let ${\cal C}$ be a function from $(u(k))_{k\in\cro{0,n+1}}$ taking its values in $`R$, and 
let $C$ and $\widehat{C}$ the two processes 
\[C(k)={\cal C}(u(k)) ~~~\textrm{ and }~~~\widehat{C}(k)={\cal C}(F_T(k))~~~\textrm{ for any }~~~k\in\cro{0,n},\]
and interpolated between integer points. Let $c_n(t)=C(nt)/k(n)$ and $\widehat{c}_n(t)={\widehat{C}(2nt)}/{k(n)}$,  
\begin{pro} Assume $({\rm H}_1)$. If $k(n)\to+\infty$, if $\dis \frac{\log n\sup_l |C(l)-C(l+1)|}{k(n)}\proba 0$,
then
\[(c_n(t))_{t\in[0,1]}\dd (c(t))_{t\in[0,1]}\]
is equivalent to
\[ ({\widehat{c}_n(t)})_{t\in[0,1]}\dd  (c(t))_{t\in[0,1]}.\]
\end{pro}
\proof  We will prove that $\|c_n-\widehat{c_n}\|_{\infty}\proba 0$. 

For any $l\in\cro{0,n}$, let $j_l$ be the integer such that $m_T(j_l)+1\leq 2l \leq m_T(j_l+1)-1$. 
Since $\widehat{C}(m_T(j_l))=C(j_l)$, we have
\be
\max_l\l|\widehat{C}(2l)-C(l)\r|\leq \max_l|\widehat{C}(2l)-\widehat{C}(m_T(j_l))|+\max_l|{C}(j_l)-C(l)|.
\ee 
Take $\rho$ large enough s.t.  $\Omega_n^\rho=\{T, \sup_l \l||u(l)|-|u(l+1)|\r|\leq \rho \log n\}$ satisfies $`P_n(\Omega_n^\rho)\to 1$ (see Lemma \ref{degdeg}). 
On $\Omega_n^\rho$, $|m_T(j_l)-2l|\leq |m_T(j_l)-m_T(j_l+1)|\leq |H(j_l)-H(j_l+1)|+2\leq  \rho\log n+2$ and then  \begin{equation}
\max_l|\widehat{C}(2l)-C\widehat(m_T(j_l))|\leq (\rho\log n+2)\sup_l |C(l)-C(l+1)|
\end{equation}
since for any $i$, $|\widehat{C}(i)-\widehat{C}(i+1)|=|{C}(j)-{C}(j+1)|$ for some $j$ (among $F_T(i)$ and $F_T(i+1)$, one is the father of the other one). Now, $|j_l-l|\leq \max_i |H_T(i)|$, and then since $\bh_n\dd \bh$, for any $`e>0$, 
\begin{equation}\label{qez}{n^{-1/2-`e}}\sup_{l}|j_l-l|\proba 0.\end{equation}
Assume that $c_n\dd c$. By \eref{qez}, 
$k(n)^{-1}\max_l|{C}(j_l)-C(l)|\proba 0,$
and then $\|c_n-\widehat{c_n}\|_{\infty}\proba 0$. \par
If $\widehat{c_n}\dd c$ then $k^{-1}(n) \max_l|\widehat{C}(2l)-\widehat{C}(m_T(j_l))|\proba 0$; on the other hand, $\max_l|C(j_l)-C(l)|=\max_l |\widehat{C}(m_T(j_l))-\widehat{C}(m_T(l))|$. Since $\max_l|m_T(j_l)-m_T(l)|\leq 2\max_l|j_l-l|+2\max_l |H_{j_l}-H_l|$, $n^{-1/2-`e}\max_l|m_T(j_l)-m_T(l)|\proba 0$. Hence
 $\max_l|C(j_l)-C(l)|/k(n)\proba 0$ and then $\|c_n-\widehat{c_n}\|_{\infty}\proba 0$. ~~~$\Box$

\subsection*{Acknowledgments} I warmly thank Svante Janson who points out some errors and imprecisions on a preliminary version of this paper and who makes numerous remarks and suggestions.
\small
\bibliographystyle{plain}

\end{document}